\documentclass[10pt,oneside,reqno,leqno,b5paper]{amsart}
\usepackage{amssymb,amsmath,latexsym}
\usepackage{amsmath,amsfonts,amsthm,latexsym}
\usepackage{mathrsfs, textcomp}
\usepackage{color,enumerate,times}
\usepackage{tikz}
\usetikzlibrary{positioning}
\usepackage{hyperref}
\usetikzlibrary{cd}
\newtheorem{theorem}{Theorem}[subsection]
\theoremstyle{definition}
\newtheorem{lem}[theorem]{Lemma}
\newtheorem{prop}[theorem]{Proposition}
\newtheorem{cor}[theorem]{Corollary}
\newtheorem{defn}[theorem]{Definition}
\newtheorem{Rem}[theorem]{Remark}
\newtheorem{exam}[theorem]{Example}
\newtheorem{thm}[theorem]{Theorem}

\newcommand{\si}{\sigma}

\usepackage{amssymb}
\usepackage{amscd}
\usepackage{eqlist, color}
\usepackage{mathrsfs, textcomp}
\usepackage{color, eqlist}
\usepackage{graphicx}
\usepackage{bm}
\usepackage{amsmath,amsfonts,amsthm,latexsym}
\usepackage[all]{xy}

\def\ifif {if and only if\ \ }
\DeclareMathOperator{\Min}{Min}

\DeclareMathOperator{\Max}{Max}
\DeclareMathOperator{\Spec}{Spec}

\DeclareMathOperator{\Ann}{Ann}
\DeclareMathOperator{\inte}{int}

\DeclareMathOperator{\cl}{cl}

\newcommand{\CO}{\mathcal{O}}

\newcommand{\sub}{\supseteq}

\newcommand{\we}{\wedge}
\newcommand{\ve}{\vee}
\newcommand{\bcap}{\bigcap}
\newcommand{\bcup}{\bigcup}

\newcommand{\ov}{\overline}

\newcommand{\LRTA}{\Leftrightarrow}

\newcommand{\dto}{\downarrow}

\newcommand{\coz}{{\rm Coz}}

\def\ifif {if and only if\ \ }
\def\sub {\subseteq}

\def\set {\setminus }

\begin{document}
\title[The class of $DA$-lattices and duality]
{Some New Classes of Lattices Related to Normal and Co-normal Lattices}

\author{Ali~Rezaei ~Aliabad}
\address{Department of Mathematics, Sahid Chamran University, Ahvaz, Iran}
\email{aliabadyr@scu.ac.ir}

\author{Navid Sabzali}
\address{Department of Mathematics, education Farhangian University, P.O. Box $14665-889$, Tehran, Iran}
\email{Navid.Sabzali@cfu.ac.ir}

\author{Ali Taherifar}
\address{Department of Mathematics, Yasouj University, Yasouj, Iran}
\email{ataherifar@yu.ac.ir, ataherifar54@gmail.com}

\begin{abstract}
We introduce two new classes of bounded lattices, called $DA$-lattices
and $DDA$-lattices, which are closely related to the classical notions
of co-normal and normal lattices, respectively. A lattice $L$ is called
a $DA$-lattice if, whenever $a\wedge b=0$, there exist
$a_1\in \bot(a)$ and $b_1\in \bot(b)$ such that
$\bot(a_1)\cap \bot(b_1)=\{0\}$, where
$\bot(x)=\{y\in L:y\wedge x=0\}$, for $x\in L$. The dual notion is called a
$DDA$-lattice. We investigate the relationships between these classes
and several other lattice-theoretic properties. In particular, we observe
that every co-normal lattice is a $DA$-lattice, whereas every normal
lattice is a $DDA$-lattice. We also study the lattice
$L^\perp=\{\bot(a):a\in L\}$ of principal annihilator ideals and prove
that $L$ is a $DA$-lattice if and only if $L^\perp$ is a $DDA$-lattice.
Moreover, we characterize quasi-complemented lattices in terms of the
corresponding $\top$-property of $L^\perp$. We show that, in general,
a lattice property and its dual need not be preserved under the
passage from $L$ to $L^\perp$. For a completely regular Hausdorff space $X$, we
prove that $Z(X)$ is a $DDA$-lattice \ifif $Z^\circ(X)$ is a $DDA$-lattice if and only if $X$ is a
$WCF$-space. Finally, for a distributive lattice $L$, we establish
that $L$ is a $DDA$-lattice if and only if $\Max L$ is Hausdorff,
equivalently, normal. Examples are provided to illustrate and delimit our results.
\end{abstract}
\maketitle
\section{Introduction}

Normality is one of the fundamental separation properties in topology and has several important algebraic counterparts in lattice theory. In particular, the notion of a normal lattice provides an algebraic framework for studying normality through distributive lattices and their spectra. The notion of a co-normal lattice is closely related to our definition
of a $DA$-lattice. Co-normal lattices were introduced by Cornish under
the name of normal lattices and have been further studied in connection
with minimal prime ideals, annihilator ideals, and topological
representation theorems. For a distributive lattice $L$, co-normality
admits several equivalent characterizations in terms of minimal prime
ideals, while normality admits corresponding characterizations in terms
of maximal prime ideals. Moreover, if $X$ is a $T_1$-space, then the
lattice $\mathcal{C}(X)$ of closed subsets of $X$ is co-normal if and only
if $X$ is normal. These classical results motivate the study of
$DA$- and $DDA$-lattices and their connections with topological
properties.

The purpose of this paper is to introduce and investigate two new classes of lattices which arise naturally from the annihilator operation. For a non-empty subset $A$ of a bounded lattice  $L$, the set \[\bot(A)=\{x\in L: x\wedge a=0, \forall a\in A\},\] is called the \emph{annihilator ideal} of $A$ in $L$. We use $\bot(a)$ (or $a^\bot$) instead of $\bot(\{a\})$. Clearly $\bot(A)=\bot({\dto}A)$ for every $A\sub L$, and if $L$ is a distributive lattice,  then for every $A\sub L$, we have $\bot(A)=\bot(<A>)$ where $<A>$ is the ideal in $L$ induced by $A$. The dual notion of $\bot(A)$ is denoted by $\top(A)$ which is the set \[\top(A)=\{x\in L: x\vee a=1, \forall a\in A\}.\] We also use $\top(a)$ (or $a^\top$) instead of $\top(\{a\})$. If $L$ is a bounded distributive lattice, then $\bot(A)$ (resp., $\top(A)$) is indeed an ideal (resp., filter) of $L$. We call a bounded lattice $L$ a $DA$-lattice if, for every $a,b\in L$ satisfying $a\wedge b=0$, there exist $a_1\in \bot(a)$ and $b_1\in \bot(b)$ such that
$\bot(a_1)\cap \bot(b_1)=\{0\}$.
The dual notion is called a $DDA$-lattice. Thus, $L$ is a $DDA$-lattice if, whenever $a\vee b=1$, there exist $a_1,b_1\in L$ such that
\[\top(a_1)\cap \top(b_1)=\{1\}
\quad\text{and}\quad
a\vee a_1=b\vee b_1=1.\]
The terminology is motivated by the fact that $DA$ and $DDA$ are dual to one another. These definitions are weaker than the corresponding co-normal and normal conditions, respectively, and therefore provide new classes naturally situated between several familiar classes of lattices.

The first part of the paper investigates the basic properties of $DA$-lattices. We show that every co-normal lattice is a $DA$-lattice. We also introduce a related class, called a $\perp$-lattice. A lattice $L$ is a $\perp$-lattice if, for every $a\in L$, there exists $b\in \bot(a)$ such that
$\bot(b)=\bot(\bot(a))=\bot^2(a)$.
Equivalently, for every $a\in L$, there exists $b\in \bot(a)$ such that
$\bot(a)\cap \bot(b)={0}$.
Every $\perp$-lattice is a $DA$-lattice, and every pseudocomplemented lattice is a $\perp$-lattice. We also obtain connections with dense and quasi-complemented lattices. In particular, for a distributive lattice, the $\perp$-property coincides with the quasi-complemented property. These results show that the $\perp$-lattice condition provides a natural intermediate notion between pseudocomplementation and the $DA$-property.

We next investigate topological realizations of the new classes. For a $T_1$-space $X$, we prove that the lattice $\mathcal{C}(X)$ of closed subsets of $X$ is a $DA$-lattice if and only if $X$ is normal (Proposition \ref{13}). Thus, the $DA$-property provides a lattice-theoretic formulation of normality for the lattice of closed sets. On the other hand, $\mathcal{C}(X)$ is a $\perp$-lattice if and only if $X$ is locally indiscrete (Proposition \ref{14}). The dual results are obtained for the lattice of open subsets. These characterizations illustrate that the new notions are closely connected with classical topological separation properties.

We then apply the theory to lattices naturally associated with completely regular Hausdorff spaces. In particular, we study the lattices $\operatorname{Coz}(X)$, $Z(X)$, and $Z^\circ(X)$ consisting of cozero-sets, zero-sets, and the interiors of zero-sets, respectively. We prove that several apparently different conditions are equivalent: $X$ is a $WCF$-space \ifif $\operatorname{Coz}(X)$ is a $DA$-lattice \ifif both lattices $Z(X)$ and $Z^\circ(X)$ are $DDA$-lattices (Theorem \ref{23}). Thus, the new lattice-theoretic properties give a unified description of the $WCF$-space condition.

The dual theory is developed in the next part of the paper. Using the filter
$\top(a)=\{x\in L:x\vee a=1\}$,
we introduce the corresponding notion of a $\top$-lattice. It is easily seen that $L$ is a $DA$-lattice if and only if its dual lattice is a $DDA$-lattice, while $L$ is a $\perp$-lattice if and only if its dual is a $\top$-lattice. Every normal lattice is a $DDA$-lattice, and every $\top$-lattice is a $DDA$-lattice. These duality results lead to further examples and clarify the position of the new classes among normal, co-normal, and related lattices.

Another important application concerns rings and their associated ideal lattices. For a commutative ring $R$, we study the lattice $Id(R)$ of ideals of $R$ using the Zariski topologies on $\Spec(R)$, $\Max(R)$, and $\Min(R)$.

For a reduced ring $R$, the lattice $Id(R)$ is a $DA$-lattice; in addition, several connections with co-normality and extremal disconnectedness are obtained. In contrast, the $DDA$-property is governed by the maximal spectrum. More precisely, we prove that
\[
Id(R)\text{ is a }DDA\text{-lattice}
\quad\Longleftrightarrow\quad
\Max(R)\text{ is a normal space, (Theorem \ref{25})}.
\]
This gives a topological characterization of the $DDA$-property for ideal lattices and provides a useful application of the new lattice concepts to commutative ring theory.

Finally, we associate with every distributive lattice $L$ the lattice
\[L^\bot=\{a^\bot:~a\in L\}=\{\bot(a):~a\in L\},\]
of principal annihilator ideals. It is proved that $L$ is a $DA$-lattice \ifif $L^\perp$ is a $DDA$-lattice (Theorem \ref{dube}). We show that $L^\perp$ is a distributive lattice and investigate its relationship with the lattice of annihilator ideals $A(L)$. In particular, $A(L)$ is always a Boolean algebra, while $L^\perp$ is closely related to the dual of $L$ through a natural quotient construction. We further establish that, for a distributive lattice $L$, $L$ is co-normal precisely when $L^\perp$ is a sublattice of the ideal lattice $I(L)$. 

The paper is organized as follows. Section 2 recalls the basic terminology and results concerning normal and co-normal lattices, annihilators, and the topological lattices used throughout the paper. Section 3 introduces $DA$-lattices and $\perp$-lattices and establishes their fundamental properties and topological applications. Section 4 develops the dual theory of $DDA$-lattices and $\top$-lattices and applies it to topological spaces and commutative rings. Section 5 studies the lattice $L^\perp$ and the lattice of annihilator ideals of a distributive lattice, obtaining further structural and duality results.

\section{Preliminaries}
We collect some terminology and results that will be needed in what follows. For known concepts that we do not define here, we refer the reader to \cite{G0, G} for lattice-theoretic concepts, to \cite{G, W} for topological spaces.
\subsection{Basic lattice  concepts}
A partially ordered set (poset) $L$ is called a \emph{lattice} if any two elements $a,b\in L$ have a least upper bound $a\vee b$ and a greatest lower bound $a\wedge b$. A lattice $L$ is called \emph{complete} if every subset $S\subseteq L$ has a least upper bound $\bigvee S$. A map $f:L\to M$ between lattices is called a \emph{lattice homomorphism} if 
$f(a\vee b)=f(a)\vee f(b)$ and $f(a\wedge b)=f(a)\wedge f(b)$ for all $a, b\in L$. A lattice $L$ is called \emph{distributive} if the identity
\[
a \wedge (b \vee c) = (a \wedge b) \vee (a \wedge c)
\]
holds, for all $a,b,c\in L$. The distributive law above and its dual are equivalent. If they exist, we denote the top and bottom elements of a lattice by $1$ and $0$, respectively, and call $L$ \emph{bounded} when both exist. Most of the distributive lattices considered in this paper are bounded. An element $a$ of a bounded  lattice $L$ is called \emph{complemented} if there exists an element $a'\in L$ such that:
\[
a\wedge a'=0
\quad\text{and}\quad
a\vee a'=1.
\]
If, in addition, $L$ is distributive, then such an element $a'$ is unique and is called the \emph{complement} of $a$. A \emph{Boolean algebra} is a bounded distributive lattice $L$ in which every $a \in L$ is complemented. The \emph{pseudocomplement} of an element $a\in L$, when it exists, is denoted by $a^*$. It is the largest element $b\in L$ such that $a\wedge b=0$. If every element of a lattice $L$ has a pseudocomplement, then $L$ is called a \emph{pseudocomplemented lattice}. For a lattice $(L, \leq)$, $L^{op}=(L, \geq)$.
\begin{prop}\label{dt}
\textit{(Principle of Duality).} Let $\mathcal{C}$ be a family of partially ordered sets closed under the formation of dual; that is, $X^{op} \in \mathcal{C}$ for every $X \in \mathcal{C}$. If a property $P$ holds for every member of $\mathcal{C}$, then its dual property $DP$ also holds for every member of $\mathcal{C}$.
\end{prop}
\subsection{The ideal lattice} 
An \emph{ideal} in a bounded  lattice $L$ is a subset $I \subseteq L$ satisfying:
\begin{enumerate}
\item $0 \in I$,
\item $a, b \in I \Rightarrow a \vee b \in I$, 
\item $b \leq a$ and $a \in I\Rightarrow b \in I$.
\end{enumerate}
 Clearly,  a non-empty subset $I$  of $L$ is an ideal  \ifif for every $a,b\in I$, we have $a\ve b\in I~\LRTA~a,b\in I$.
For a distributive lattice $L$, the family $I(L)$ of all ideals, ordered by inclusion, forms a frame. Meets in $I(L)$ are given by set-theoretic intersection, and the join of a family $\mathcal{J} \subseteq I(L)$ is given by the formula
\[
\bigvee\mathcal{J} = \bigcup\{I_1\vee I_2\vee \cdots\vee I_n \mid I_1,I_2,\ldots,I_n \in \mathcal{J},\, n \in \mathbb{N}\},
\]
where 
\[
I_1 \vee I_2 \vee \cdots \vee I_n = \{i_1 \vee i_2 \vee \cdots \vee i_n \mid i_k \in I_k \text{ for } 1 \leq k \leq n\}.
\]
In particular, if $\mathcal{J}$ is directed, then $\bigvee\mathcal{J} = \bigcup\mathcal{J}.$
In general, it is customary to refer to $I(L)$ as the \emph{ideal lattice} of $L$. We have the following result so.
\begin{cor}\cite[Corollary 104]{G0}
A lattice $L$ is distributive \ifif for any two ideals $I, J$ of $L$:
\[I\vee J=\{i\vee j: i\in I, j\in J\}.\]
\end{cor}

It is well-known that for $a$ in a distributive lattice $L$, $\bot(a)=\{x\in L: x\wedge a=0\}$ is an ideal. Indeed, it is an annihilator ideal. Thus, for $a, b \in L$, we write $\bot(a) \wedge \bot(b)$ for the intersection, $\bot(a) \cap \bot(b)$, of these two ideals, and $\bot(a) \vee \bot(b)$ for their join in the ideal lattice of $L$. Also, for $a\in L$, we have $\top(a)=\{x\in L: x\vee a=1\}$. It is easily seen that if $L$ is a distributive lattice, then $\top(a)$ is a filter on $L$. 

 For subsets $R, S$ of a lattice $L$, $R\overline{\wedge} S=\{a\wedge b: a\in R~\text{and}~b\in S\}$, $R\underline{\vee} S=\{a\vee b: a\in R~\text{and}~b\in S\}$ and  we know that for a subset $R$ of $L$, $\bot(R)=\{x\in L: x\wedge y=0,~ \forall y\in R\}$. It is easy to see that  $\bot(R)=\bot^3(R)$.
\begin{lem}\label{03}
Let $L$ be a lattice and $a, b\in L$.
\begin{enumerate}
\item $a\wedge b=0$ if and only if $\bot^2(a)\wedge \bot^2(b)=\{0\}$. 

\item $a\vee b=1$ \ifif $\top^2(a)\wedge \top^2(b)=\{1\}$. 

\item If $L$ is a distributive lattice, then $\bot(a\vee b)=\bot(a)\wedge \bot(b)$. 

\item If $L$ is a pseudocomplemented lattice, then $\bot(a)\wedge \bot(b)=\{0\}$ if and only if $a^*\wedge b^*=0$. 

\item A lattice $L$ is pseudocomplemented \ifif for each $a\in L$ there exists $b\in L$ such that $\bot(a)=\downarrow b$.
\end{enumerate}
\end{lem}
\begin{proof}
(1) We prove the statement in a slightly more general form. For $R,S\subseteq L$, we have
$$R\overline{\wedge} S=\{0\}~\Leftrightarrow~ R\subseteq \bot(S)\Leftrightarrow~ R\overline{\wedge} \bot^2(S)=\{0\}~\Leftrightarrow~\bot^2(S)\subseteq \bot(R)\Leftrightarrow~\bot^2(R)\overline{\wedge} \bot^2(S)=\{0\}.$$Taking $R=\{a\}$ and $S=\{b\}$ yields the desired result.

(2) The proof is evident, by the part (1) and Proposition \ref{dt}.

(3) This is immediate from distributivity.

(4)$\Rightarrow$ Assume that $c=a^*\wedge b^*$. Thus, $c\in \bot(a)\cap \bot(b)=\{0\}$ and so $c=0$. 

$\Leftarrow$ Suppose that $x\in \bot(a)\wedge \bot(b)$. Then, $x\wedge a=x\wedge b=0$. This implies $x\leq a^*\wedge b^*=0$. 

(5)$\Rightarrow$ Let $a\in L$. Set $b=a^*$. Then, we have $b\in \bot(a)$ and $x\in \bot(a)$ implies $x\leq a^*=b$. Hence, $\bot(a)=\downarrow b$.

(5)$\Leftarrow$ Let $a\in L$. By hypothesis, there exists $b\in L$ such that $\bot(a)=\downarrow b$. Clearly, $b=\max(\bot(a))=a^*$.
\end{proof}
\begin{prop}\label{242}
Let $B$ be a Boolean algebra, and let $L\subseteq B$ be a lattice with
respect to the order inherited from $B$. Put
\[
L'=\{a'\in B:a\in L\},
\]
where $a'$ denotes the complement of $a$ in $B$. Then the following
statements hold.
\begin{enumerate}
\item The ordered sets $(L,\geq)$ and $(L',\leq)$ are order-isomorphic.

\item The set $L'$ is a lattice with the operations
\[
a'\vee_{L'}b'=(a\wedge_L b)',
\qquad
a'\wedge_{L'}b'=(a\vee_L b)',
\]
for all $a,b\in L$. Moreover, $L'$ is isomorphic to the dual lattice
$L^{op}$.

\item If $L$ is a sublattice of $B$, then the lattice operations in
$L'$ coincide with those inherited from $B$; that is,
\[
a'\vee_{L'}b'=a'\vee_B b'
\qquad\text{and}\qquad
a'\wedge_{L'}b'=a'\wedge_B b'
\]
for all $a,b\in L$.
\end{enumerate}
\end{prop}

\begin{proof}
(1) Define
\[
f:L\longrightarrow L',\qquad\text{by}~\qquad f(a)=a'.
\]
Clearly, $f$ is bijective. Moreover, for all $a,b\in L$,
\[
a\leq_{L^{op}} b
\iff
b\leq_L a
\iff
a'\leq_B b'
\iff
f(a)\leq f(b).
\]
Hence, $f$ is an order isomorphism from $L^{op}$ onto $L'$.

(2) Since $f:L^{op}\to L'$ is an order isomorphism, it preserves
joins and meets. Therefore, for all $a,b\in L$,
\[
a'\vee_{L'}b'
=
f(a)\vee_{L'}f(b)
=
f(a\wedge_L b)
=
(a\wedge_L b)',
\]
and
\[
a'\wedge_{L'}b'
=
f(a)\wedge_{L'}f(b)
=
f(a\vee_L b)
=
(a\vee_L b)'.
\]
Thus $L'$ is a lattice isomorphic to $L^{op}$.

(3) Suppose that $L$ is a sublattice of $B$. Then
\[
a\wedge_L b=a\wedge_B b
\qquad\text{and}\qquad
a\vee_L b=a\vee_B b.
\]
Hence, by De Morgan's laws,
\[
a'\vee_{L'}b'
=
(a\wedge_L b)'
=
(a\wedge_B b)'
=
a'\vee_B b',
\]
and similarly,
\[
a'\wedge_{L'}b'
=
(a\vee_L b)'
=
(a\vee_B b)'
=
a'\wedge_B b'.
\]
Therefore, the operations of $L'$ agree with the operations inherited
from $B$.
\end{proof}
\subsection{Normal (resp, co-normal) lattices} Recall from \cite[3.6]{J}, a distributive lattice $L$ is called a \emph{normal lattice} if  $a, b\in L$ satisfy $a\vee b=1$,  there exist $a_1, b_1\in L$ such that $a_1\wedge b_1=0$ and $a\vee a_1=b\vee b_1=1$.

\emph{Co-normal} lattices were introduced by Cornish in \cite{C}, under the name normal lattice and characterised in \cite[Theorem 2.4 $(b)\Leftrightarrow (h)$]{C} as those distributive lattices $(L, \wedge, \vee, 0, 1)$  satisfying the condition that, for all $a, b\in L$ with $a\wedge b=0$ there exist $x, y\in L$ such that $x\vee y=1$ and $x\wedge a=y\wedge b=0$. Such lattices have been further studied in \cite{AT,BB1, MT, MT1, TA,T}. As observed in \cite{AMT}, it is preferable to call this class of lattices \emph{co-normal lattices}, since the term `normal'' refers to the dual property. Trivially, every Boolean algebra is a co-normal lattice. It was mention in \cite{C} that the lattice of closed sets of a $T_1$-space is co-normal if and only if the space is normal. We recall below some well-known results concerning normal and co-normal lattices that will be used throughout the paper.
\begin{lem}\cite[Proposition 3.7]{J}\label{simo1}
For a distributive lattice $L$ the following statements are equivalent.
\begin{enumerate}
\item $L$ is a normal lattice.

\item Every prime ideal in $L$ is contained in a unique maximal ideal.

\item The inclusion map from $\Max L$ to $\Spec L$ is a split monomorphism.
\end{enumerate}
\end{lem}
\begin{thm}\cite[3.6]{J}\label{simo2}
The following statements hold.
\begin{enumerate}
\item A topological space $X$ is normal \ifif the lattice $\Omega(X)$ ($\mathcal{O}(X)$) is normal. 

\item If $L$ is a distributive lattice, then $\Max L$ is a compact $T_1$-space. 

\item If $L$ is a normal lattice, then $\Max L$ is a Hausdorff space.
\end{enumerate}
\end{thm}
\begin{thm}\cite[Theorem 2.4]{C}\label{simo3}
The following statements for a distributive lattice $L$ are equivalent.
\begin{enumerate}
\item The lattice $L$ is co-normal.

\item Any two distant minimal prime ideals are comaximal.

\item For any $x, y\in L$, $x\wedge y=0$ implies $\bot(x)$ and $\bot(y)$ are comaximal.

\item For any $x, y\in L$, $\bot(x\wedge y)=\bot(x)\vee \bot(y)$.

\item Each maximal ideal contains a unique minimal prime ideal. 

\item Each prime ideal contains a unique minimal prime ideal.
\end{enumerate}
\end{thm}
\subsection{Rings of Continuous Functions and Topological Concepts}
In this paper, $C(X)$ (resp., $C^{*}(X)$) is denoted by the ring of all (resp., all bounded) real-valued continuous functions on a completely regular Hausdorff space $X$. For each $f \in C(X)$, the set $f^{-1}(\{0\})$ is called the \emph{zero-set} of $f$, and is denoted by $Z(f)$, while  $\mathrm{Coz}(f)=X\setminus Z(f)$ is the \emph{cozero-set} of $f$. 
The set of all zero-sets in $X$, $Z(X)$, partially ordered by inclusion and with the following operations is a distributive lattice,
\[Z(f)\vee Z(g)=Z(fg),~\text{and}~Z(f)\wedge Z(g)=Z(f)\cap Z(g),~\text{where}~Z(f), Z(g)\in Z(X).\]
Also, the set of all cozero-sets in $X$, denoted by $\operatorname{Coz}(X)$, partially ordered by inclusion and equipped with the following operations, forms a distributive lattice:
\[
\operatorname{Coz}(f) \wedge \operatorname{Coz}(g) = \operatorname{Coz}(fg)
\quad \text{and} \quad
\operatorname{Coz}(f) \vee \operatorname{Coz}(g) = \operatorname{Coz}(f) \cup \operatorname{Coz}(g),~
\text{where}~f, g \in C(X).\]
 For $A\sub X$, the closure of $A$ in $X$ ($\cl_XA$) is denoted by $\ov{A}$ and the interior of $A$ in $X$ ($\inte_XA$) is denoted by $A^\circ$, see \cite{GJ}.

A space $X$ is \textit{cozero complemented} if, given any cozero
set $U$, there is a disjoint cozero set $V$ such that $U\cup V$ is
dense in $X$. There are many equivalent algebraic or topological
conditions for  cozero complemented spaces, see  \cite{HW} and \cite{T1}.

Recall from \cite{AT} that a topological space $X$ is a $WCF$-space if every pair of disjoint cozero-sets in $X$ can be separated by two disjoint $Z^\circ$-sets. That is, if $\mathrm{Coz}(f)$ and $\mathrm{Coz}(g)$ are two disjoint cozero-sets in $X$, then there exist two zero-sets $Z(f_1)$ and $Z(f_2)$ containing $\mathrm{Coz}(f)$ and $\mathrm{Coz}(g)$, respectively such that $Z^\circ(f_1)\cap Z^\circ(f_2)=\emptyset$.

Recall from \cite{G} that if  $L$ is a pseudocomplement lattice and $a,b\in L$, then $(a\wedge b)^{**}=a^{**}\wedge b^{**}$. We apply this fact in the following result to the lattice of open (resp., closed) subsets of a topological space $X$, denoted by $\mathcal{O}(X)$ (resp., $\mathcal{C}(X)$). Before, we note that if $X$ is a topological space and $U\in\CO(X)$, then $U^{**}=\inte_X\cl_X(U)$. The lattice of open subsets is also denoted by $\Omega(X)$.
\begin{lem}\label{00}
The following statements hold.
\begin{enumerate}
\item If $A, B\in \CO(X)$ (resp., $A, B\in \mathcal{C}(X)$), then \[\inte_X\cl_X(A\cap B)=\inte_X\cl_X A\cap \inte_X\cl_X B ~(\text{resp}., \cl_X\inte_X(A\cup B)=\cl_X\inte_X A\cup\cl_X\inte_X B).\]

\item $U$ and $V$ are  disjoint open sets in a space $X$, \ifif $\inte\cl(U)\cap\inte\cl(V)=\emptyset$. 
\end{enumerate}
\end{lem}
\begin{lem}\label{razagh}
Let $X$ be a topological space. Suppose $\mathcal{U} \subseteq \mathcal{O}(X)$ (respectively, $\mathcal{F} \subseteq \mathcal{C}(X)$) is a sublattice of $\mathcal{O}(X)$ (respectively, of $\mathcal{C}(X)$). Then the following statements hold.
\begin{enumerate}
\item The set $\overline{\mathcal{U}} = \{ \overline{U} : U \in \mathcal{U} \}$, equipped with the operations
\[
\overline{U} \wedge \overline{V} = \overline{U \cap V}
\quad \text{and} \quad
\overline{U} \vee \overline{V} = \overline{U \cup V}
\]
for all $U, V \in \mathcal{U}$, is a bounded distributive lattice.

\item The set $\mathcal{F}^{\circ} = \{ F^{\circ} : F \in \mathcal{F} \}$, equipped with the operations
\[
F^{\circ} \wedge H^{\circ} = (F \cap H)^{\circ}
\quad \text{and} \quad
F^{\circ} \vee H^{\circ} = (F \cup H)^{\circ}
\]
for all $F, H \in \mathcal{F}$, is a bounded distributive lattice.

\item The set $\overline{\mathcal{U}}^{^\circ} = \{ \overline{U}^{^\circ} : U \in \mathcal{U} \}$, equipped with the operations
\[
\overline{U}^{^\circ} \wedge \overline{V}^{^\circ} = \overline{U \cap V}^{^\circ}
\quad \text{and} \quad
\overline{U}^{^\circ} \vee \overline{V}^{^\circ} = \overline{U \cup V}^{^\circ}
\]
for all $U, V \in \mathcal{U}$, is a bounded distributive lattice. Indeed, $\overline{\mathcal{U}}^{^\circ}$ is a sublattice of the lattice of regular open sets.

\item The set $\overline{\mathcal{F}^{^\circ}} = \{ \overline{F^{^\circ}} : F \in \mathcal{F} \}$, equipped with the operations
\[
\overline{F^{^\circ}} \wedge \overline{H^{^\circ}} = \overline{(F \cap H)^{^\circ}}
\quad \text{and} \quad
\overline{F^{^\circ}} \vee \overline{H^{^\circ}} = \overline{(F \cup H)^{^\circ}}
\]
for all $F, H \in \mathcal{F}$, is a bounded distributive lattice. Indeed, $\overline{\mathcal{F}^{^\circ}}$ is a sublattice of the lattice of regular closed sets.
\end{enumerate}
\end{lem}

\begin{proof}
(1) Clearly, $\overline{U \cap V}$ is a lower bound of $\{ \overline{U}, \overline{V} \}$. Now let $\overline{W}$ be a lower bound of $\{ \overline{U}, \overline{V} \}$, where $W \in \mathcal{U}$. Then we have
\[
\overline{W} \subseteq \overline{U} \cap \overline{V}.
\]
Taking interiors and closures, we obtain
\[
W \subseteq \overline{W}^{^\circ} \subseteq \overline{U}^{^\circ} \cap \overline{V}^{^\circ}
= \overline{U \cap V}^{^\circ} \subseteq \overline{U \cap V}.
\]
Thus $\overline{W} \subseteq \overline{U \cap V}$, and hence $\overline{U \cap V} = \overline{U} \wedge \overline{V}$.

Since $\overline{U \cup V} = \overline{U} \cup \overline{V}$, the second operation is immediate. For distributivity, let $U, V, W \in \mathcal{U}$. Then
\[
\overline{U} \wedge (\overline{V} \vee \overline{W})
= \overline{U} \wedge \overline{V \cup W}
= \overline{U \cap (V \cup W)}
= \overline{(U \cap V) \cup (U \cap W)}
= \overline{U \cap V} \cup \overline{U \cap W}
= (\overline{U} \wedge \overline{V}) \vee (\overline{U} \wedge \overline{W}).
\]
Hence $\overline{\mathcal{U}}$ is distributive. The bounds are $\overline{\varnothing} = \varnothing$ and $\overline{X} = X$, so it is bounded.

(2) The proof is similar to (1).

(3) The proof is similar to (1).

(4) The proof is similar to (1).
\end{proof}
\begin{lem}\label{02}
Let $X$ be a topological space. Then, the following statements hold:
\begin{enumerate}
\item The set $\overline{\operatorname{Coz}(X)} = \{ \overline{\operatorname{Coz}(f)} : f \in C(X) \}$, partially ordered by inclusion, is a bounded distributive lattice with the following operations:
\[
\overline{\operatorname{Coz}(f)} \vee \overline{\operatorname{Coz}(g)} = \overline{\operatorname{Coz}(f^2 + g^2)}
\quad \text{and} \quad
\overline{\operatorname{Coz}(f)} \wedge \overline{\operatorname{Coz}(g)} = \overline{\operatorname{Coz}(fg)},
\]
for all $f, g \in C(X)$.

\item The set $Z^\circ(X) = \{ Z^\circ(f) : f \in C(X) \}$, partially ordered by inclusion, is a bounded distributive lattice with the following operations:
\[
Z^\circ(f) \wedge Z^\circ(g) = Z^\circ(f^2 + g^2)
\quad \text{and} \quad
Z^\circ(f) \vee Z^\circ(g) = Z^\circ(fg),
\]
for all $f, g \in C(X)$.
\end{enumerate}
\end{lem}

\begin{proof}
(1) This follows immediately from Part (1) of Lemma~\ref{razagh}.

(2) This follows from Part (1) by taking complements. See also \cite{TA}.
\end{proof}

\section{$DA$-lattice}
\subsection{DA-lattice and $\perp$-lattice}
\begin{defn}
A Lattice $L$  is called a $DA$-lattice if for each $a, b\in L$ with $a\wedge b=0$, there is  $(a_1, b_1)\in \bot(a)\times \bot(b)$ such that $\bot(a_1)\wedge \bot(b_1)=\{0\}$.
\end{defn}
\begin{defn} 
A Lattice $L$ is called a $\perp$-lattice if for every $a\in L$, there exists $b\in \bot(a)$ such that \[\bot(b)=\bot(\bot(a))=\bot^2(a).\]
\end{defn}
Clearly, a lattice $L$ is a $\perp$-lattice \ifif for each $a\in L$ there exists $b\in \bot(a)$ such that $\bot(a)\wedge \bot(b)=\{0\}$.
\begin{lem}\label{11}
The following statements hold.
\begin{enumerate}
\item Every co-normal lattice is a $DA$-lattice.

\item If $L$ is a $\perp$-lattice, then it is a $DA$-lattice.

\item If $L$ is a pseudocomplemented lattice, then $L$ is a $\perp$-lattice and so a $DA$-lattice.

\item Every dense lattice is a $\perp$-lattice and hence a $DA$-lattice.

\item For every distributive lattice $L$, the ideal lattice, $I(L)$, is a pseudocomplemented lattice and so a $DA$-lattice.
\end{enumerate}
\end{lem}
\begin{proof}
(1)
Assume $a, b\in L$ with $a\wedge b=0$. Then, by hypothesis, there are $a_1, b_1\in L$ such that $a_1\vee b_1=1$ and  $a\wedge a_1=b\wedge b_1=0$. Since $L$ is a distributive lattice, $a_1\vee b_1=1$ implies  $\bot(a_1)\wedge \bot(b_1)=\bot(a_1\vee b_1)=0$, by Lemma \ref{03}(2). Hence, $L$ is a $DA$-lattice.

(2) Suppose that $L$ is a $\perp$-lattice and $a\wedge b=0$. By the hypothesis, there exist $(a_1,b_1)\in \bot(a)\times \bot(b)$ such that $\bot(a_1)=\bot^2(a)$ and $\bot(b_1)=\bot^2(b)$. Clearly, $a\in \bot(a_1)$, $b\in \bot(b_1)$ and we have  \[\bot(a_1)\cap \bot(b_1)=\bot^2(a)\cap \bot^2(b)=\{0\},\] by Lemma \ref{03}(1). Hence, $L$ is a $DA$-lattice.

(3) By Lemma \ref{03},  $\bot(a^*)=\bot({\dto}a^*)=\bot(\bot(a))=\bot^2(a)$, so we are done.

(4) Assume that $L$ is a dense lattice and $x\in L$. If $x=0$, then $1\in \bot(0)$ and $\bot(1)=\{0\}=\bot^2(0)$, and so we are done. Now, assume that $x\neq 0$. Then, by hypothesis, $\bot(x)=\{0\}$ and $\bot(0)=\bot^2(x)$.

(5) By \cite[Lemma 3.2.1(1)]{T}, for each $I\in I(L)$, we have $I^*=\bot(I)$. Thus, $I(L)$ is a pseudo-complemented lattice, see also \cite[Lemma 1.2]{G}. So, we are done by Part (3).
\end{proof}
\begin{exam}\label{12}
(1) It is well-known that $Z(X)$, the lattice of zero-sets in a space X, is a co-normal lattice, see Lemma 1.4 in \cite{BB1}. Hence, by Lemma \ref{11}(1), $Z(X)$ is a $DA$-lattice.  This provides an example of a distributive co-normal lattice that is a $DA$-lattice but is not a $\perp$-lattice,  see Corollary \ref{cor17}.
 
(2) Let $X$ be an infinite set equipped with the cofinite topology. Then the proper closed subsets of $X$ are precisely the finite subsets of $X$. Consider the lattice
\[
\mathcal{C}(X)=\{F\subseteq X:F\text{ is finite}\}\cup\{X\}.
\]
For every closed subset $A$ of $X$, we have
\[
\bot(A)=\{F\in\mathcal{C}(X):F\cap A=\emptyset\}.
\]
In particular, for closed subsets $A, B$ of $X$,
\[
\bot(A)\cap \bot(B)=\bot(A\cup B).
\]

We claim that $\mathcal{C}(X)$ is not a DA-lattice. Indeed, let $A,B\in\mathcal{C}(X)$ be nonempty disjoint finite subsets of $X$. Suppose that $A_1\in \bot(A)$ and $B_1\in \bot(B)$. Since $A_1$ and $B_1$ are finite, their union $A_1\cup B_1$ is finite. Because $X$ is infinite, there exists a nonempty finite subset $H\subseteq X$ such that
$
H\cap(A_1\cup B_1)=\emptyset.
$
Consequently,
$
H\in \bot(A_1)\cap \bot(B_1),
$
and hence
$
\bot(A_1)\cap \bot(B_1)\neq\{0\}.
$
Therefore, $\mathcal{C}(X)$ does not satisfy the defining condition of a $DA$-lattice.

(3) Consider the following lattice L: 

\[\begin{tikzpicture}[
    node distance = 0.8cm,
    every node/.style = {circle, draw, minimum size=0.2cm, inner sep=2pt, fill=white},
    edge/.style = {draw, thick}
]

\node (1)   at (0, 3)   {1};
\node (c)   at (0, 2)   {c};
\node (a)   at (-1.5, 1) {a};
\node (b)   at (1.5, 1)  {b};
\node (0)   at (0, 0)   {0};

\draw[edge] (0) -- (a);
\draw[edge] (0) -- (b);
\draw[edge] (a) -- (c);
\draw[edge] (b) -- (c);
\draw[edge] (c) -- (1);

\end{tikzpicture}\]
Then, it is easily verified  that $L$ is not a co-normal lattice. On the other hand, we have $\bot(a)=\{b, 0\}$ and $\bot(b)=\{a, 0\}$. Thus, if we take $a_1=b$ and $b_1=a$, then we have $\bot(a_1)\wedge \bot(b_1)=\{0\}$, $a\in \bot(a_1)$ and $b\in \bot(b_1)$. Hence, $L$ is a $DA$-lattice. In fact, $L$ is a $\perp$-lattice.
\end{exam}
 Since, for every $a, b\in L$, we have $\bot(a\vee b)\subseteq \bot(a)\wedge \bot(b)$, it follows that every $\perp$-lattice is a quasi-complemented, see \cite[Definition 5.1]{C}.
Moreover, if $L$ is a distributive lattice, then, $L$ is a $\perp$-lattice \ifif $L$ is a quasi-complemented lattice. However, the following example shows that these two notions are distinct in general.
\begin{exam}
Let $L$ be the following lattice. 
\[\begin{tikzpicture}[
    node distance = 0.8cm,
    every node/.style = {circle, draw, minimum size=0.2cm, inner sep=2pt, fill=white},
    edge/.style = {draw, thick}
]

\node (1)   at (0, 3.5)  {1};
\node (a)   at (-2, 1.8) {a};
\node (b)   at (2, 1.8)  {b};
\node (c)   at (0, 0.8)  {c};
\node (0)   at (0, 0)    {0};

\draw[edge] (0) -- (a);
\draw[edge] (0) -- (b);
\draw[edge] (0) -- (c);

\draw[edge] (a) -- (1);
\draw[edge] (b) -- (1);
\draw[edge] (c) -- (1);

\end{tikzpicture}\]
Then, it is easy to see that $L$ is a quasi-complemented lattice. But, it is not a $\perp$-lattice. In fact, it is not a $DA$-lattice. For, we have $\bot(a)\wedge \bot(b)\neq 0$,  for each $a, b\in L$. 
\end{exam}
\begin{prop}\label{13}
For a $T_1$-space $X$, the following statements are equivalent.
\begin{enumerate}
\item The lattice $\mathcal{C}(X)$ is a $DA$-lattice.

\item $X$ is a normal space.

\item The lattice $\mathcal{C}(X)$ is a co-normal lattice.
\end{enumerate}
\end{prop}
\begin{proof}
(1)$\Rightarrow$(2) Let $A$ and $B$ be two disjoint closed subsets of $X$. By hypothesis, there are two closed subsets $A_1, B_1$ in $X$ such that $\bot(A_1\cup B_1)=\bot(A_1)\wedge \bot(B_1)=0$, $A\in \bot(A_1)$ and $B\in \bot(B_1)$. We claim that $A_1\cup B_1=X$ and so we are done. To see it, let $A_1\cup B_1$ be a proper subset of $X$. Then, there exists $x\not\in A_1\cup B_1$. Thus, $x\not\in A_1$ and $x\not\in B_1$. That is $\{x\}\in \bot(A_1)\cap \bot(B_1)$, a contradiction.

(2)$\Rightarrow$(3)  See \cite[Proposition 4.4]{C}.

(3)$\Rightarrow$(1) This follows from Lemma \ref{11}.
\end{proof}
Recall that a \emph{locally indiscrete} space is a topological space in which every open set is also a closed set. A discrete space is obviously locally indiscrete. A typical example of a nontrivial locally indiscrete space is $\mathrm{Min}(R)$ with the Zariski topology, where $R$ is a principal ideal ring, and, of course, $\mathrm{Min}(R)$ is the set of all minimal prime ideals of the ring $R$.
\begin{prop}\label{14}
The following statements hold.
\begin{enumerate}
\item For every space $X$, the lattice $\CO(X)$ ($\Omega(X)$) is a  pseudocomplemented lattice and hence a $\perp$-lattice.

\item For a $T_1$-space, the  lattice $\mathcal{C}(X)$ is a $\perp$-lattice \ifif $X$ is a locally indiscrete space.

\end{enumerate}
\end{prop}
\begin{proof}
(1) Let $U$ be an open subset of $X$. Then, it is easy to see that $U^*=X\setminus\cl_X U=\overline{U}^{^c}\in\CO(X)$. Hence, $\CO(X)$ is a pseudocomplemented lattice.

(2) $\Rightarrow$  Let $A\in\mathcal{C}(X)$. By hypothesis, there exists a closed set $B$ such that $B\in \bot(A)$ (i.e., $B\cap A=\emptyset$) and $\bot(A\cup B)=\bot(A) \cap \bot(B)=\{0\}$. The last equality implies $A\cup B=X$. We also have $A\cap B=\emptyset$, hence $A=X\setminus B$ is open, i.e., $X$ is a locally indiscrete.

$\Leftarrow$ Let $A\in\mathcal{C}(X)$. By hypothesis, $A$ is open and hence $A^c$ is closed. We have $A^c\in \bot(A)$ and $\bot(A)\wedge \bot(A^c)=\bot(A\cup A^c)=\{0\}$. Hence, $\mathcal{C}(X)$ is a $\perp$-lattice.
\end{proof}
Consider a normal space which is not locally indiscrete (e.g., $\Bbb R$ with the standard topology). Propositions \ref{13} and \ref{14} show that $\mathcal{C}(X)$ is a co-normal lattice but is not a $\perp$-lattice.
\begin{lem}\label{aa}
Let $X$ be a topological space. Suppose $\mathcal{U} \subseteq \mathcal{O}(X)$ is a sublattice of $\mathcal{O}(X)$ and a base for the open sets of $X$ (respectively, $\mathcal{F} \subseteq \mathcal{C}(X)$ is a sublattice of $\mathcal{C}(X)$ such that $\mathcal{F}$ is a nhood-base for the open sets). Then the following statements hold:
\begin{enumerate}
\item For every $U \in \mathcal{U}$, we have $\bot(U) = \{0\}$ in $\mathcal{U}$ if and only if $\overline{U} = X$.

\item For every $F \in \mathcal{F}$, we have $\bot(F) = \{0\}$ in $\mathcal{F}$ if and only if $F = X$.

\item For every $U \in \mathcal{U}$, we have $\bot(\overline{U}) = \{0\}$ in $\overline{\mathcal{U}}$ if and only if $\overline{U} = X$.

\item For every $F \in \mathcal{F}$, we have $\bot(F^{\circ}) = \{0\}$ in $\mathcal{F}^{\circ}$ if and only if $\overline{F^{\circ}} = X$.
\end{enumerate}
\end{lem}

\begin{proof}
(1) $(\Rightarrow)$ Suppose $U \in \mathcal{U}$  such that $\bot(U) = \{0\}$. Let $\emptyset \neq V \in \mathcal{O}(X)$. Since $\mathcal{U}$ is a base for the open sets, there exists $\emptyset \neq W \in \mathcal{U}$ such that $W \subseteq V$. As $W \notin \bot(U)$, we have $W \cap U = W \wedge U \neq \emptyset$. Hence every non-empty open set $V$ meets $U$, so $\overline{U} = X$.

$(\Leftarrow)$ Suppose $U \in \mathcal{U}$  such that $\overline{U} = X$. Then $V \cap U \neq \emptyset$ for every non-empty open set $V$. In particular, for every non-empty $W \in \mathcal{U}$, we have $W \cap U \neq \emptyset$. Hence $\bot(U) = \{0\}$ in $\mathcal{U}$.

(2) $(\Rightarrow)$ Suppose $F \in \mathcal{F}$ with $\bot(F) = \{0\}$. Assume, for a contradiction, that $X \setminus F \neq \emptyset$. Since $X \setminus F$ is a non-empty open set and $\mathcal{F}$ is a nhood base for the open sets, there exists $ \emptyset \neq H\in \mathcal{F}$ such that  $H\subseteq X\setminus F$. Thus, $H \cap F = \emptyset$, so $\emptyset\neq H \in \bot(F)$, contradicting $\bot(F) = \{0\}$. Therefore $F = X$.

$(\Leftarrow)$ If $F = X$, then clearly $\bot(F) = \{0\}$ in $\mathcal{F}$.

(3) $(\Rightarrow)$ Suppose $U \in \mathcal{U}$  such that $\bot(\overline{U}) = \{0\}$ in $\overline{\mathcal{U}}$. Let $ \emptyset\neq V \in \mathcal{O}(X)$. Since $\mathcal{U}$ is a base, there exists $\emptyset \neq W \in \mathcal{U}$ such that $W \subseteq V$. As $\overline{W} \notin \bot(\overline{U})$, we have $\overline{W} \wedge \overline{U} = \overline{W \cap U} \neq \emptyset$. Thus $\emptyset \neq W \cap U \subseteq V \cap U$. Hence every non-empty open set $V$ meets $U$, so $\overline{U} = X$.

$(\Leftarrow)$ Suppose $U \in \mathcal{U}$  such that $\overline{U} = X$. Let $\overline{V} \in \overline{\mathcal{U}}$ be non-empty, where $V \in \mathcal{U}$. Since $\overline{U} = X$ and $V$ is non-empty open, we have $U \cap V \neq \emptyset$, so $\overline{V} \wedge \overline{U} = \overline{V \cap U} \neq \emptyset$. Hence $\bot(\overline{U}) = \{0\}$ in $\overline{\mathcal{U}}$.

(4) $(\Rightarrow)$ Suppose $F \in \mathcal{F}$  such that $\bot(F^{\circ}) = \{0\}$ in $\mathcal{F}^{\circ}$. Let $\emptyset\neq V \in \mathcal{O}(X)$. Since $\mathcal{F}$ is a nhood base for open sets, there exists $H \in \mathcal{F}$ such that for $x\in V$, $x\in H^\circ\subseteq H \subseteq V$. As $H^\circ \notin \bot(F^{\circ})$, we have $H^\circ \cap F^{\circ} = H^\circ \wedge F^{\circ} \neq \emptyset$. Thus $H^\circ \cap F^{\circ}$ is a non-empty subset of $V \cap F^{\circ}$. Hence every non-empty open set meets $F^{\circ}$, so $\overline{F^{\circ}} = X$.

$(\Leftarrow)$ Suppose $F \in \mathcal{F}$  such that $\overline{F^{\circ}} = X$. Let $H \in \mathcal{F}$ be such that $H^{\circ} \neq \emptyset$. Since $\overline{F^{\circ}} = X$ and $H^{\circ}$ is a non-empty open set, we have $H^{\circ} \cap F^{\circ} \neq \emptyset$. Thus
\[
H^{\circ} \wedge F^{\circ} = (H \cap F)^{\circ} = H^{\circ} \cap F^{\circ} \neq \emptyset.
\]
Therefore $\bot(F^{\circ}) = \{0\}$ in $\mathcal{F}^{\circ}$.
\end{proof}
\begin{prop}\label{fati}
Let $X$ be a topological space. Suppose $\mathcal{U} \subseteq \mathcal{O}(X)$ is a sublattice of $\mathcal{O}(X)$ and a base for the open sets of $X$ (respectively, $\mathcal{F} \subseteq \mathcal{C}(X)$ is a sublattice of $\mathcal{C}(X)$ such that $\mathcal{F}$ is a nhood base for the open sets). Then the following statements hold:
\begin{enumerate}
\item $\mathcal{U}$ is a $\perp$-lattice if and only if for every $U \in \mathcal{U}$, there exists $V \in \mathcal{U}$ such that $U \cap V = \emptyset$ and $\overline{U \cup V} = X$.

\item $\overline{\mathcal{U}}$ is a $\perp$-lattice if and only if  $\overline{\mathcal{U}}$ is a Boolean lattice.

\item $\mathcal{F}$ is a $\perp$-lattice if and only if $\mathcal{F}$ is Boolean.

\item $\mathcal{F}^{\circ}$ is a $\perp$-lattice if and only if for every $F \in \mathcal{F}$, there exists $H \in \mathcal{F}$ such that $F^{\circ} \cap H^{\circ} = \emptyset$ and $\overline{(F \cup H)^{\circ}} = X$.
\end{enumerate}
\end{prop}

\begin{proof}
(1) By Lemma \ref{aa}, we have the following equivalences:
\[
\begin{aligned}
\mathcal{U} \text{ is a } \perp\text{-lattice}
&\Longleftrightarrow \forall U \in \mathcal{U}, \; \exists V \in \bot(U) \text{ such that } \bot(U \cup V) = \bot(U) \cap \bot(V) = \{0\} \\
&\Longleftrightarrow \forall U \in \mathcal{U}, \; \exists V \in \mathcal{U} \text{ such that } U \cap V = U \wedge V = \emptyset \text{ and } \overline{U \cup V} = X.
\end{aligned}
\]
The last equality follows from Lemma~\ref{aa}.

(2) By Lemma~\ref{aa}, \[
\begin{aligned}
\ov{\mathcal{U}} \text{ is a } \perp\text{-lattice}
&\Longleftrightarrow \forall U \in \mathcal{U}, \; \exists \ov{V} \in \bot(\ov{U}) \text{ such that }~ \bot(\ov{U\cup V})=\bot(\ov{U} \cup \ov{V}) = \bot(U) \cap \bot(V) = \{0\} \\
&\Longleftrightarrow \forall U \in \mathcal{U}, \; \exists V \in \mathcal{U} \text{ such that } \ov{U} \cap \ov{V} = \ov{U} \wedge \ov{V} = \emptyset \text{ and } \overline{U \cup V} = X.
\end{aligned}
\] That is each $\ov{U}\in\ov{\mathcal{U}}$ has complemented.  Hence, $\ov{\mathcal{U}}$ is a Boolean algebra.

(3) By Lemma \ref{aa}, we have:
\[
\begin{aligned}
\mathcal{F} \text{ is a } \perp\text{-lattice}
&\Longleftrightarrow \forall F \in \mathcal{F}, \; \exists H \in \bot(F) \text{ such that } \bot(F \cup H) = \bot(F) \wedge \bot(H) = \{0\} \\
&\Longleftrightarrow \forall F \in \mathcal{F}, \; \exists H \in \mathcal{F} \text{ such that } F \cap H = F \wedge H = \emptyset \text{ and } F \cup H = X.
\end{aligned}
\]
That is, every element of $\mathcal{F}$ is complemented. Hence $\mathcal{F}$ is Boolean.

(4) By Lemma \ref{aa}, we have:
\[
\begin{aligned}
\mathcal{F}^{\circ} \text{ is a } \perp\text{-lattice}
&\Longleftrightarrow \forall F \in \mathcal{F}, \; \exists H \in \mathcal{F} \text{ with } H^{\circ} \in \bot(F^{\circ}) \text{ such that } \bot((F \cup H)^{\circ}) = \bot(F^{\circ}) \cap \bot(H^{\circ}) = \{0\} \\
&\Longleftrightarrow \forall F \in \mathcal{F}, \; \exists H \in \mathcal{F} \text{ such that } F^{\circ} \cap H^{\circ} = F^{\circ} \wedge H^{\circ} = \emptyset \text{ and } \overline{(F \cup H)^{\circ}} = X.
\end{aligned}
\]
This completes the proof.
\end{proof}
\begin{cor}\cite[Proposition 3.7]{S}\label{cor17}
The lattice $Z(X)$ is a $\perp$-lattice if and only if $X$ is $P$-space.
\end{cor}
\begin{proof}
This follows from  Proposition \ref{fati}(3) and this fact that $X$ is a $P$-space \ifif every zero-set is open, see \cite[4J]{GJ}.
\end{proof}
\begin{prop}\label{sadri}
For a completely regular Hausdorff space $X$, the following statements are equivalent.
\begin{enumerate}
\item The lattice $\coz(X)$ is a $\perp$-lattice.

\item The space $X$ is a cozero-complemented space.

\item The lattice $Z^\circ(X)$ is a $\perp$-lattice.

\item The lattice $\ov{\coz(X)}$ is a $\perp$-lattice.
\end{enumerate}
\end{prop}
\begin{proof}
(1)$\Leftrightarrow$(2) This follows from Proposition \ref{fati}(1).

(2)$\Leftrightarrow$(3) This follows from Proposition \ref{fati}(4).

(2)$\Leftrightarrow$(4) This follows from Proposition \ref{fati}(2).
\end{proof}

\section{Duality for $DA$-lattices and $\perp$-lattices}
Recall that for $a\in L$,  $\top(a)=\{x\in L: x\vee a=1\}$. It is easy to see that $\top(a)$ is a filter on $L$. For each $a\in L$, we have $\top(\bot(a))=\{1\}~\text{and}~ \bot(\top(a))=\{0\}$.

\subsection{DDA lattice, $\top$-lattice and their applications}
\begin{defn}
A lattice $L$ is called a \emph{dual $DA$-lattice} (briefly, a \emph{$DDA$-lattice}) if for each $a, b \in L$ with $a \vee b = 1$, there exist elements $a_1, b_1 \in L$ such that
$\top(a_1) \wedge \top(b_1) = \{1\}$
and $a \vee a_1 = b \vee b_1 = 1$.
\end{defn}
\begin{defn}
A lattice $L$ is called a $\top$-lattice if for each $a\in L$, there exists  $b\in \top(a)$ such that $\top^2(a)=\top(b)$.
\end{defn}
By Lemma \ref{03}, it is easy to see that a lattice $L$ is a $\top$-lattice \ifif for each $a\in L$, there exists $b\in \top(a)$ such that $\top(a)\cap \top(b)=\{1\}$. Trivially, $L$ is a $DA$-lattice \ifif $L^{op}$ is a $DDA$-lattice. Also, $L$ is $\perp$-lattice \ifif $L^{op}$ is a $\top$-lattice. Using these facts and Lemma \ref{11} to obtain the next results.
\begin{lem}
The following statements hold.
\begin{enumerate}
\item Every normal lattice is a $DDA$-lattice.  

\item  Every $\top$-lattice is a $DDA$-lattice.

\item Every co-dense lattice is a $\top$-lattice and hence a $DDA$-lattice.

\item Every Boolean algebra is a $\top$-lattice and hence a $DDA$-lattice.
\end{enumerate}
\end{lem}
\begin{exam}\label{21}
(1) Since the lattice $Z(X)$ is co-normal, by Proposition \ref{242},
the lattice $\coz(X)$ is normal and hence is a $DDA$-lattice.

(2) Using the dual of item (3) of Example \ref{12}, we obtain a
distributive lattice $L$ that is not normal, whereas $L$ is a
$\top$-lattice and hence a $DDA$-lattice.
\end{exam} 
\begin{prop}
For a $T_1$-space $X$, the following  statements hold.
\begin{enumerate}
\item $\mathcal{O}(X)$ is a $DDA$-lattice.

\item $X$ is a normal lattice.

\item $\mathcal{O}(X)$ is a normal lattice.
\end{enumerate}
\end{prop}
The following result follows from Propositions \ref{14} and \ref{sadri}.
\begin{prop}\label{22}
The following statements hold.
\begin{enumerate}
\item  For a $T_1$-space $X$, the lattice $\CO(X)$ is a $\top$-lattice \ifif $X$ is a locally indiscrete space.

\item For any   space $X$, the lattice $\mathcal{C}(X)$ is a $\top$-lattice and hence is a $DDA$-lattice.

\item For a completely regular space $X$, the lattice $Z(X)$ is a $\top$-lattice \ifif $X$ is a cozero-complemented lattice. 
\end{enumerate}
\end{prop}
We now apply the lattice properties $DA$ and $DDA$ to obtain a characterization of the class of $WCF$-spaces.
\begin{thm}\label{23}
For a completely regular Hausdorff space $X$, the following statements are equivalent.
\begin{enumerate}
\item The space $X$ is a $WCF$-space.

\item The lattice $\coz(X)$ is a $DA$-lattice.

\item The lattice $\ov{\coz(X)}$ is a $DA$-lattice.

\item The lattice $Z(X)$ is a $DDA$-lattice.

\item The lattice $Z^{\circ}(X)$ is a $DDA$-lattice.

\end{enumerate}
\end{thm}
\begin{proof}
(1)$\Rightarrow$(2) Let $A, B\in\coz(X)$ with $A\cap B=\emptyset$. Then, by hypothesis, there exist $Z_1, Z_2\in Z(X)$ such that 
\[A\subseteq Z_1,~ B\subseteq Z_2~\text{and}~ Z_1^{\circ}\cap Z_2^{\circ}=\emptyset.\] The last equality says that $\ov{X\set(Z_1\cap Z_2)}=X$. 
Set $U=X\set(Z_1\cap Z_2)$, by Lemma \ref{aa},  we have,
\[\bot((X\setminus Z_1)\vee (X\setminus Z_2))=\bot(U)=\emptyset.\] 
Hence, $\coz(X)$ is a $DA$-lattice.

(2)$\Rightarrow$(3) Let $\ov{\coz(f)}, \ov{\coz(g)}\in\ov{\coz(X)}$ such that $\ov{\coz(fg)}=\ov{\coz(f)}\we \ov{\coz(g)}=\emptyset$. Then, we have 
\[\coz(f)\we\coz(g)=\coz(f)\cap\coz(g)=\coz(fg)=\emptyset.\] 
By hypothesis, there are $\coz(f_1), \coz(g_1)\in\coz(X)$ such that 
\[\bot(\coz(f_1)\cup \coz(g_1))=\{\emptyset\}~\text{and}~ \coz(ff_1)=\coz(f)\we \coz(f_1)=\coz(g)\we \coz(g_1)=\coz(gg_1)=\emptyset.\] 
This together with Lemma \ref{aa}, implies
\[\ov{\coz(f_1)}\cup \ov{\coz(g_1)}=(\ov{\coz(f_1)\cup\coz(g_1)}=X~~,~~ \ov{\coz(f)}\we \ov{\coz(f_1)}=\ov{\coz(ff_1)}=\emptyset~~,\text{and}\]
\[\ov{\coz(g)}\we \ov{\coz(g_1)}=\ov{\coz(gg_1)}=\emptyset.\] 
Therefore, $\ov{\coz(X)}$ is a $DA$-lattice. 

(3)$\Rightarrow$(1) Let $\coz(f)$ and $\coz(g)$ be two disjoint cozero-sets in $X$. Then, \[\ov{\coz(f)}\we\ov{\coz(g)}=\coz(f)\cap\coz(g)=\ov{\coz(fg)}=\emptyset.\]
Since $\ov{\coz(X)}$ is a $DA$-lattice,  there exist $\ov{\coz(f_1)}, \ov{\coz(g_1)}\in\ov{\coz(X)}$ such that 
\[\bot(\ov{\coz(f_1)}\vee\ov{\coz(g_1)})=\bot(\ov{\coz(f_1)}\cup\ov{\coz(g_1)})=\{\emptyset\}~\text{and}~\ov{\coz(ff_1)}=\]
\[\ov{\coz(f)}\we\ov{\coz(f_1)}=\ov{\coz(g)}\we\ov{\coz(g_1)}=\ov{\coz(gg_1)}=\emptyset.\] 
By Lemma \ref{aa}, it follows that $\ov{\coz(f_1)}\cup\ov{\coz(g_1)}=X$ and $\coz(ff_1)=\coz(gg_1)=\emptyset$. Hence, $\ov{\coz(f_1)}\cup\ov{\coz(g_1)}=X$ and $\coz(f)\cap \coz(f_1)=\coz(g)\cap \coz(g_1)=\emptyset$. This implies 
\[\coz(f)\subseteq Z(f_1),~\coz(g)\subseteq Z(g_1)~ \text{and}~Z^{\circ}(f_1)\cap Z^{\circ}(g_1)=\emptyset.\] Therefore, $X$ is a $WCF$-space.

(2)$\Leftrightarrow$(4) It is easy to see that the lattices $(Z(X), \subseteq)$ and $\coz(X)^{op}$ are isomorphic and we have the lattice  $\coz(X)$ is $DA$ \ifif the lattice $\coz(X)^{op}$ is    $DDA$, so we are done. 

(3)$\Leftrightarrow$(5) It is easy to see that the lattices $(Z^{\circ}(X), \subseteq)$ and $\ov{\coz(X)}^{op}$ are isomorphic and we have $\ov{\coz(X)}$ is $DA$ \ifif the lattice $\ov{\coz(X)}^{op}$ is $DDA$, this completes the proof.

\end{proof}
\begin{exam}
(1) Let $X$ be a non-$WCF$-space (e.g., see \cite[Example 3.9]{AT} ). Then, $Z(X)$ is not a $DDA$-lattice, by Theorem \ref{23}.

(2) It is proved  in \cite{ST} that $Z(X)$ is a normal lattice \ifif $X$ is an $F$-space. Now, consider a $WCF$-space which is not an $F$-space (e.g., \cite[Example 3.9]{AT}). Then, $Z(X)$ is a $DDA$-space which is not a normal lattice.
\end{exam}
\subsection{Lattice properties and Zariski topology}
In this subsection all rings are commutative and with identity. We use $\Spec(R), \Max(R)~  \text{and}~ \Min(R)$ for the
spaces of prime ideals, maximal ideals and  minimal prime ideals of $R$, respectively with the
$hull$-$kernel$ topology. For a subset $S$ of $R$, let $h(S)=\{P\in\Spec(R): S\subseteq P\}$. If $S=\{a\}$, then we use $h(a)$. The set $\{h(a): a\in R\}$ forms a base for closed sets in $\Spec(R)$. $\Max(R)$ and $\Min(R)$ are subspaces of $\Spec(R)$, and we use $h_M(S)$ and $h_{m}(S)$ instead of  $h(S)\cap\Max(R)$ and $h(S)\cap\Min(R)$, respectively. The lattice of ideals of a ring $R$ is denoted by $Id(R)$. A reduced ring $R$ is a ring in which for $a\in R$, $a^2=0$ implies $a=0$. Evidently, $R$ is a reduced ring \ifif for each $I\in Id(R)$, $I\cap \Ann(I)=0$. 
Whenever $R$ is a reduced ring, $\operatorname{Id}(R)$ is a pseudocomplemented lattice and, consequently, a $\perp$-lattice. Indeed, the reducedness hypothesis implies that for each $I \in \operatorname{Id}(R)$, we have $I \cap \operatorname{Ann}(I) = 0$. Moreover, for any $J \in \operatorname{Id}(R)$ with $J \cap I = 0$, we have $J \subseteq \operatorname{Ann}(I)$; hence $I^* = \operatorname{Ann}(I)$. This shows that $\operatorname{Id}(R)$ is a pseudocomplemented lattice. However, the converse is not true. For instance, consider $R = \mathbb{Z}_4$. Clearly, $\operatorname{Id}(R)$ is a pseudocomplemented lattice, whereas $R$ is not reduced. 
\begin{thm}\label{19}
Let $R$ be a reduced ring. Then, the following statements are equivalent.
\begin{enumerate}
\item The lattice $Id(R)$ is a co-normal lattice. 

\item The space $\Min(R)$ is extremally disconnected.

\item The space $\Spec(R)$ is extremally disconnected.
\end{enumerate}
\end{thm}
\begin{proof}
(1)$\Rightarrow$(2) Let $U$ and $V$ be two disjoint open sets in $\Min(R)$. Then, there are two ideals $I, J$ of $R$ such that $U=h_m^{c}(I)$ and $V=h_m^{c}(J)$ and $h_m^{c}(I)\cap h_m^{c}(J)=\emptyset$. This implies $h_m^{c}(IJ)=\emptyset$. Since $R$ is a reduced ring, it follows that $I\cap J=IJ=0$. By hypothesis, there are two ideals $I_1, J_1$ of $R$ such that \[I_1\vee J_1=I_1+ J_1=R~\text{and}~I\cap I_1=J\cap J_1=0.\] This implies \[h_m(I_1)\cap h_m(J_1)=h_m(I_1+J_1)=\emptyset~\text{and}~h_m^{c}(I)\cap h_m^{c}(I_1)=h_m^{c}(J)\cap h_m^{c}(J_1)=\emptyset.\] Hence, 
\[U=h_m^{c}(I)\subseteq h_m(I_1)~\text{and}~V=h_m^{c}(J)\subseteq h_m(J_1)~\text{and}~h_m(I_1)\cap h_m(J_1)=h_m(I_1+J_1)=\emptyset.\] Therefore, $\cl(U)\cap \cl(V)=\emptyset$. By \cite[1H]{GJ}, $\Min(R)$ is an extremally disconnected space. 

(2)$\Rightarrow$(3) Let $I, J\in Id(R)$ with $I\cap J=0$.  Then, \[h_m^{c}(I)\cap h_m^{c}(J)=h_m^{c}(I\cap J)=\emptyset.\] By hypothesis and \cite[1H]{GJ}, we have $\cl(h_m^{c}(I))\cap\cl(h_m^{c}(J)=\emptyset$. On the other hand, there are two ideals $I_1, J_1$ of $R$ such that \[\cl(h_m^{c}(I))=h_m(I_1)~\text{and}~ \cl(h_m^{c}(J))=h_m(J_1).\] Thus, we have \[h_m(I_1+ J_1)=h_m(I_1)\cap h_m(J_1)=\emptyset,\]and\[h_m(I\cap I_1)=h_m(I_1)\cup h_m(I)=\Min(R)=h_m(J)\cup h_m(J_1)=h_m(J\cap J_1).\] This implies $I_1+ J_1=R$ and $I\cap I_1=J\cap J_1=0$. Hence, \[I_1\subseteq\Ann(I),~J_1\subseteq\Ann(J)~ \text{and}~I_1+ J_1=R.\] Therefore, $\Ann(I)+ \Ann(J)=R$. By \cite[Theorem 2.14]{ATT}, $\Spec(R)$ is extremally disconnected. 

(3)$\Rightarrow$(1)  Let $I, J\in Id(R)$ with $I\cap J=0$. By hypothesis and \cite[Theorem 2.14]{ATT}, $\Ann(I)+ \Ann(J)=R$. Put $I_1=\Ann(I)$ and $J_1=\Ann(J)$. So, we are done. 
\end{proof}
By \cite[Theorem 2.14]{ATT}, $C(X)$ is a Baer ring \ifif $\Spec(C(X))$ is extremally disconnected space. It also is well-known that $C(X)$ is Baer ring \ifif $X$ is an extremally disconnected space (Theorem 2.12 in \cite{T1}). Hence, by Theorem \ref{19}, $Id(C(X))$ is a co-normal lattice \ifif $X$ is an extremally disconnected space. On the other hand, we have seen that $Id(R)$ for any reduced ring $R$ is a $DA$-lattice. If we consider a non-extremally disconnected space $X$ (e.g., the real number with the standard topology), then $Id(C(X))$ is a $DA$-lattice which is not a co-normal lattice.
We need the following easy lemma in the sequel.
\begin{lem}\label{24a}
For $I\in Id(R)$, $\top(I)=\{R\}$ in the lattice $Id(R)$ \ifif $I\subseteq J(R)$, where $J(R)=\bcap\Max(R)$ is the radical Jacobson  of the ring $R$.
\end{lem}
\begin{thm}\label{25}
For a ring $R$, the lattice $Id(R)$ is a $DDA$-lattice if and only if the space $\Max(R)$ is a normal ($T_4$) space.
\end{thm}

\begin{proof}
$(\Rightarrow)$
It is enough to show that $\Max(R)$ is Hausdorff, since $\Max(R)$ is compact. Let $M_1,M_2\in\Max(R)$ with $M_1\neq M_2$. Then
$M_1\vee M_2=M_1+M_2=R$.
Since $Id(R)$ is a $DDA$-lattice, there exist $I_1,J_1\in Id(R)$ such that
\[
\top(I_1\cap J_1)=\top(I_1)\cap\top(J_1)=\{R\}
\]
and
\[
M_1+I_1=M_2+J_1=R.
\]
By Lemma \ref{24a},
\[
\top(I_1\cap J_1)=\{R\}
\quad\Longrightarrow\quad
I_1\cap J_1\subseteq J(R).
\]
Consequently,
\[
h_M(I_1\cap J_1)=\Max(R).
\]
Since
\[
h_M(I_1\cap J_1)=h_M(I_1)\cup h_M(J_1),
\]
we obtain
\[
h_M^c(I_1)\cap h_M^c(J_1)=\emptyset.
\]
Moreover, $M_1+I_1=R$ and $M_2+J_1=R$ imply
\[
M_1\notin h_M(I_1)
\qquad\text{and}\qquad
M_2\notin h_M(J_1).
\]
Hence
\[
M_1\in h_M^c(I_1),\qquad
M_2\in h_M^c(J_1),
\]
where $h_M^c(I_1)$ and $h_M^c(J_1)$ are disjoint open neighborhoods of $M_1$ and $M_2$, respectively. Thus $\Max(R)$ is Hausdorff. Since $\Max(R)$ is compact, it is normal.

$(\Leftarrow)$
Let $I,J\in Id(R)$ satisfy
$I+J=I\vee J=R$.
Then,
\[
h_M(I)\cap h_M(J)
=h_M(I+J)
=h_M(R)
=\emptyset.
\]
Thus $h_M(I)$ and $h_M(J)$ are disjoint closed subsets of $\Max(R)$. Since $\Max(R)$ is normal, there exist disjoint open subsets $U,V$ of $\Max(R)$ such that
\[
h_M(I)\subseteq U,\qquad\text{and}\qquad h_M(J)\subseteq V.
\]
By the definition of the topology on $\Max(R)$, there exist $I_1,J_1\in Id(R)$ such that
\[
U=h_M^c(I_1),\qquad
V=h_M^c(J_1).
\]
Therefore,
\[
h_M(I)\subseteq h_M^c(I_1),
\qquad
h_M(J)\subseteq h_M^c(J_1).
\]
Since $U\cap V=\emptyset$, we have
\[
h_M^c(I_1)\cap h_M^c(J_1)=\emptyset.
\]
Hence
\[
h_M(I_1\cap J_1)
=h_M(I_1)\cup h_M(J_1)
=\Max(R),
\]
and consequently
\[
I_1\cap J_1\subseteq J(R).
\]
By Lemma \ref{24a}, it follows that
\[
\top(I_1\cap J_1)=\{R\}.
\]

Furthermore,
\[
h_M(I+I_1)
=h_M(I)\cap h_M(I_1)
=\emptyset,
\]
because $h_M(I)\subseteq h_M^c(I_1)$. Thus
$
I+I_1=R$.
Similarly,
\[
h_M(J+J_1)
=h_M(J)\cap h_M(J_1)
=\emptyset,
\]
and hence
$J+J_1=R$.
Therefore, for every $I,J\in Id(R)$ with $I\vee J=R$, there exist $I_1,J_1\in Id(R)$ such that
\[\top(I_1)\cap\top(J_1)=\{R\}~\text{and}
~I\vee I_1=J\vee J_1=R.
\]
Thus $Id(R)$ is a $DDA$-lattice.
\end{proof}
Our next goal is to characterize the $DDA$-lattice property for a distributive lattice $L$. We first recall the spectrum topology. For a distributive lattice $L$, the set $\operatorname{Pt}(I(L))$ of points of $I(L)$ consists of the proper prime ideals of $L$; that is, ideals $I \subsetneq L$ such that $x \wedge y \in I$ implies $x \in I$ or $y \in I$. The \emph{spectrum} of $L$, denoted by $\Spec L$, is the space whose underlying set is the collection of all proper prime ideals of $L$, equipped with the \emph{hull-kernel} topology --- also called the \emph{spectral} or \emph{Jacobson} topology --- in which the open sets are
\[
h^c(I) = \{ P \in \Spec L : P \not\supseteq I \},
\]
(see \cite{AL} and \cite{J}). In fact, the family $\{ h(a) : a \in L \}$, where $h(a)=\{P\in\Spec L: a\in P\}$, forms a base for the closed sets of this space. The space $\Max L$ is a subspace of $\Spec L$, and its closed base elements are of the form $h_M(a) = h(a) \cap \Max L$. 

If $L$ is a normal lattice, then $\Max L$ is a Hausdorff space, by Theorem~\ref{simo3}. However, the converse is not true in general. We will see that this fact is equivalent to $L$ being a $DDA$-lattice. First, we need the following lemma.
\begin{lem}\label{26}
For $a\in L$, we have $\top(a)=\{1\}$ \ifif $h_M^c(a)=\emptyset$ in $\Max L$.
\end{lem}
\begin{proof}
($\Rightarrow$). If $h_M^c(a)\neq\emptyset$, then there exists $M\in\Max L$ such that $a\not\in M$. Thus, the ideal generated by $M\cup\{a\}$ properly contains $M$ and so is whole of the $L$. Hence, there exists $b\in M$ such that $1=a\vee b$. That is $b\neq 1\in \top(a)$, a contradiction.

($\Leftarrow$). Suppose that $\top(a)\neq\{1\}$. Thus, there exists $b\neq 1\in L$ such that $a\ve b=1$. Let $M\in\Max(L)$ be such $b\in M$. Clearly, $a\notin M$ and so $h_M^c(a)\neq\emptyset$.
\end{proof}
\begin{thm}\label{27}
A bounded distributive lattice $L$ is a $DDA$-lattice if and only if $\Max L$ is a Hausdorff space (equivalently, a normal space).
\end{thm}

\begin{proof}
$(\Rightarrow)$
Let $M_1,M_2\in\Max L$ with $M_1\neq M_2$. Choose
$a\in M_1\setminus M_2$. Since $M_2$ is a maximal ideal, the ideal
generated by $M_2\cup\{a\}$ is equal to $L$. Hence, there exists
$b\in M_2$ such that
\[
a\vee b=1.
\]
Since $L$ is a $DDA$-lattice, there exist $a_1,b_1\in L$ such that
\[
\top(a_1\wedge b_1)=\top(a_1)\cap\top(b_1)=\{1\}~\text{and}~a\vee a_1=b\vee b_1=1.
\]
By Lemma \ref{26},
\[
h_M^c(a_1\wedge b_1)=\emptyset.
\]
Since
\[
h_M^c(a_1\wedge b_1)
=h_M^c(a_1)\cap h_M^c(b_1),
\]
we obtain
\[
h_M^c(a_1)\cap h_M^c(b_1)=\emptyset.
\]
Furthermore, since $a\vee a_1=1$ and $b\vee b_1=1$,
\[
h_M(a)\subseteq h_M^c(a_1)
\qquad\text{and}\qquad
h_M(b)\subseteq h_M^c(b_1).
\]
Since $a\in M_1$ and $b\in M_2$, we have
\[
M_1\in h_M(a)\subseteq h_M^c(a_1)~\text{and}~
M_2\in h_M(b)\subseteq h_M^c(b_1).
\]
Thus, $h_M^c(a_1)$ and $h_M^c(b_1)$ are disjoint open
neighborhoods of $M_1$ and $M_2$, respectively. Therefore, $\Max L$
is Hausdorff. Since $\Max L$ is a compact $T_1$-space, it follows that $\Max L$ is
normal.

$(\Leftarrow)$
Suppose that $\Max L$ is normal. Let $a,b\in L$ satisfy $a\vee b=1$.
Then
\[
h_M(a)\cap h_M(b)=h_M(a\vee b)=h_M(1)=\emptyset.
\]
Thus, $h_M(a)$ and $h_M(b)$ are disjoint closed subsets of $\Max L$.
By normality, there exist disjoint open subsets $U,V$ of $\Max L$
such that
\[
h_M(a)\subseteq U
\qquad\text{and}\qquad
h_M(b)\subseteq V.
\]
Hence, there exist ideals $I,J$ of $L$ such that
\[
U=h_M^c(I),\qquad V=h_M^c(J).
\]
Consequently,
\[
h_M(a)\subseteq h_M^c(I),
\qquad
h_M(b)\subseteq h_M^c(J),
\]
and
\[
h_M^c(I)\cap h_M^c(J)=\emptyset.
\]
Therefore,
\[
h_M^c(I\cap J)
=h_M^c(I)\cap h_M^c(J)
=\emptyset,
\]
and so $I\cap J=0$.
By Lemma \ref{26}, it follows that
$
I\cap J=0$.

Since $\Max L$ is compact and $h_M(a)\subseteq h_M^c(I)=\bcup_{x\in I}h_M^c(x)$, there exist
$a_1,\ldots,a_n\in I$ such that
\[
h_M(a)\subseteq
h_M^c(a_1)\cup\cdots\cup h_M^c(a_n)
=h_M^c(a_1\wedge\cdots\wedge a_n).
\]
Put
\[
a_0=a_1\wedge a_2\wedge...\wedge a_n\in I.
\]
Then
\[
h_M(a)\subseteq h_M^c(a_0).
\]
Similarly, there exists $b_0\in J$ such that
\[
h_M(b)\subseteq h_M^c(b_0).
\]
Thus,
\[
h_M(a\vee a_0)
=h_M(a)\cap h_M(a_0)
=\emptyset,
\]
and hence $a\vee a_0=1$.
Similarly, $b\vee b_0=1$.
Moreover, since $a_0\in I$, $b_0\in J$, and $I\cap J=\{0\}$, we have $a_0\wedge b_0=0$.
Therefore,
\[
\top(a_0\wedge b_0)=1,~\text{and}~
a\vee a_0=b\vee b_0=1.
\]
Hence $L$ is a $DDA$-lattice.
\end{proof}

\section{A New lattice from a lattice L}
Let $L$ be a bounded lattice. In this section, especially in the first subsection, we focus our main attention on the lattice $L^\bot=\{a^\bot:~a\in L\}=\{\bot(a):~a\in L\}$.
\subsection{The lattice $L^\perp$}
\begin{prop}
Let $L$ be a lattice and let $L^\perp = \{ \bot(a) : a \in L \}$. Then the following statements hold:
\begin{enumerate}
\item $L^\perp$, with the natural order, is a lattice with the following operations:
\[
\bot(a) \vee_{L^\perp} \bot(b) = \bot(a \wedge b)
\quad \text{and} \quad
\bot(a) \wedge_{L^\perp} \bot(b) = \bot(a \vee b)
\]
for all $a, b \in L$.

\item If $L$ is distributive, then for all $a, b \in L$,
\[
\bot(a) \vee_{L^\perp} \bot(b) = \bot(a \wedge b)
\quad \text{and} \quad
\bot(a) \wedge_{L^\perp} \bot(b) = \bot(a) \cap \bot(b) = \bot(a \vee b).
\]

\item If $L$ is distributive, then $L^\perp$ is a distributive lattice.
\end{enumerate}
\end{prop}

\begin{proof}
(1) We have
\[
\bot(a) \subseteq \bot(a \wedge b) \quad \text{and} \quad \bot(b) \subseteq \bot(a \wedge b).
\]
Now suppose $\bot(a) \subseteq \bot(x)$ and $\bot(b) \subseteq \bot(x)$ for some $x \in L$. We claim that $\bot(a \wedge b) \subseteq \bot(x)$. To prove this, let $r \in \bot(a \wedge b)$. Then $r \wedge (a \wedge b) = 0$. Hence, $(r \wedge b) \wedge a = 0$, so $r \wedge b \in \bot(a) \subseteq \bot(x)$. So, $r \wedge b \wedge x = 0$. This implies $(r \wedge x) \wedge b = 0$, so $r \wedge x \in \bot(b) \subseteq \bot(x)$. Therefore, $r \wedge x = (r \wedge x) \wedge x = 0$, which means $r \in \bot(x)$. Thus $\bot(a \wedge b) \subseteq \bot(x)$, proving that $\bot(a) \vee_{L^\perp} \bot(b) = \bot(a \wedge b)$. The proof for the meet operation is similar.

(2) It remains only to show that $\bot(a) \cap \bot(b) = \bot(a \vee b)$. The inclusion $\bot(a \vee b) \subseteq \bot(a) \cap \bot(b)$ is clear. Conversely, suppose $x \in \bot(a) \cap \bot(b)$. Then $x \wedge a = 0$ and $x \wedge b = 0$. By distributivity,
\[
x \wedge (a \vee b) = (x \wedge a) \vee (x \wedge b) = 0 \vee 0 = 0.
\]
Thus $x \in \bot(a \vee b)$. Hence $\bot(a) \cap \bot(b) \subseteq \bot(a \vee b)$, and so equality holds.

(3) Assume $\bot(a), \bot(b), \bot(c) \in L^\perp$. Then, using the operations established in (1) and (2), we have
\[
\begin{aligned}
\bot(a) \wedge_{L^\perp} (\bot(b) \vee_{L^\perp} \bot(c))
&= \bot(a) \wedge_{L^\perp} \bot(b \wedge c) \\
&= \bot(a \vee (b \wedge c)) \\
&= \bot((a \vee b) \wedge (a \vee c)) \\
&= \bot(a \vee b) \vee_{L^\perp} \bot(a \vee c) \\
&= (\bot(a) \wedge_{L^\perp} \bot(b)) \vee_{L^\perp} (\bot(a) \wedge_{L^\perp} \bot(c)).
\end{aligned}
\]
Hence $L^\perp$ is distributive.
\end{proof}

\begin{Rem}
Note that if $L^\perp$ is distributive, then $L$ is not necessarily distributive. For example, let $L = N_5 = \{ 0, a, b, c, 1 \}$, where $0 < a < b < 1$, $0 < c < 1$ and $a$ and $b$ are not comparable to $c$. Then
\[
L^\perp = \{ \bot(0), \bot(c), \bot(b), \bot(1) \}
= \{ \{0\}, \{0,c\}, \{0,a,b,1\}, \{0,a,b,c,1\} \} \cong M_2,
\]
which is distributive, whereas $L$ is not.
\end{Rem}
\begin{prop}\label{sadr}
For a distributive lattice $L$, the following statements are equivalent.
\begin{enumerate}
\item The lattice $L^\perp$ is a sublattice of $I(L)$

\item For each $a, b\in L$, we have $\bot(a)\vee \bot(b)=\bot(a\wedge b)$.

\item The lattice $L$ is a co-normal lattice.
\end{enumerate}
\end{prop}
\begin{proof}
For equivalency of (2) and (3) see Theorem \ref{simo3}.  Trivially (2) implies (1). To prove that (1) implies (2), it is enough to show that $\bot(a\wedge b)\subseteq \bot(a)\vee \bot(b)$. Since,   the inequality $\bot(a)\vee \bot(b) \subseteq \bot(a\wedge b)=\bot(a) \vee_{L^\perp} \bot(b)$ always holds. By hypothesis, $\bot(a)\vee \bot(b)\in L^\perp$ and it is containing both $\bot(a)$ and $\bot(b)$. Thus, \[\bot(a) \vee_{L^\perp} \bot(b)=\bot(a\wedge b)\subseteq \bot(a)\vee \bot(b),\] so we are done.
\end{proof}
\begin{thm}\label{dube}
Let $L$ be a distributive lattice. 
\begin{enumerate}
\item If $L$ is a co-normal lattice, then  $L^\perp$ is a normal lattice.

\item $\top(\bot(a))=\{\bot(b): b\in \bot(a)\}$ in $L^\perp$, for every $a\in L$.

\item $L$ is a $DA$-lattice \ifif $L^\perp$ is a $DDA$-lattice.

\item $L$ is a quasi-complemented lattice ($\perp$-lattice) \ifif $L^\perp$ is a $\top$-lattice.

\end{enumerate}
\end{thm}
\begin{proof}
(1) Let $L$ be a co-normal lattice and let $\bot(a), \bot(b) \in L^\perp$ with $\bot(a) \vee_{L^\perp} \bot(b) = 1_{L^\perp} = L$. Then $\bot(a \wedge b) = 1_{L^\perp} = L$, which implies $a \wedge b = 0$. By hypothesis, there exist $a_1, b_1 \in L$ such that $a_1 \vee b_1 = 1$ and $a \wedge a_1 = b \wedge b_1 = 0$. Thus we have
\[
\bot(a_1) \wedge \bot(b_1) = \bot(a_1 \vee b_1) = \bot(1) = \{0\},
\]
and
\[
\bot(a) \vee_{L^\perp} \bot(a_1) = \bot(a \wedge a_1) = \bot(0) = L = 1_{L^\perp},
\]
and similarly,
\[
\bot(b) \vee_{L^\perp} \bot(b_1) = \bot(b \wedge b_1) = \bot(0) = L = 1_{L^\perp}.
\]
Hence $L^\perp$ is a normal lattice.

(2) Suppose $a \in L$. Then clearly
\[
\top(\bot(a)) = \{ \bot(b) : L = \bot(b) \vee_{L^\perp} \bot(a) = \bot(b \wedge a) \}
= \{ \bot(b) : b \wedge a = 0 \}
= \{ \bot(b) : b \in \bot(a) \}.
\]

(3) First, assume $L$ is a $DA$-lattice and $\bot(a), \bot(b) \in L^\perp$ with $\bot(a) \vee_{L^\perp} \bot(b) = L = 1_{L^\perp}$. Then $\bot(a \wedge b) = L$, and so $a \wedge b = 0$. By hypothesis, there exist $a_1, b_1 \in L$ such that
\[
\bot(a_1) \wedge \bot(b_1) = \bot(a_1 \vee b_1) = \{0\}
\]
and
\[
a \wedge a_1 = b \wedge b_1 = 0.
\]
Thus we have
\[
\bot(a) \vee_{L^\perp} \bot(a_1) = \bot(a \wedge a_1) = \bot(0) = L
\]
and
\[
\bot(b) \vee_{L^\perp} \bot(b_1) = \bot(b \wedge b_1) = \bot(0) = L.
\]
Now let $\bot(x) \in \top(\bot(a_1)) \cap \top(\bot(b_1))$. Then, by part (2), $x \in \bot(a_1) \cap \bot(b_1) = \{0\}$, and so $\bot(x) = L$. This implies
\[
\top(\bot(a_1)) \cap \top(\bot(b_1)) = L = 1_{L^\perp}.
\]

Conversely, let $a \wedge b = 0$ in $L$. Then
\[
\bot(a) \vee_{L^\perp} \bot(b) = \bot(a \wedge b) = \bot(0) = L = 1_{L^\perp}.
\]
By hypothesis, there exist $\bot(a_1), \bot(b_1) \in L^\perp$ such that
\[
\top(\bot(a_1)) \cap \top(\bot(b_1)) = \{1_{L^\perp}\},
\]
\[
\bot(a \wedge a_1) = \bot(a) \vee_{L^\perp} \bot(a_1) = L = 1_{L^\perp},
\]
and
\[
\bot(b \wedge b_1) = \bot(b) \vee_{L^\perp} \bot(b_1) = L = 1_{L^\perp}.
\]
Thus $a \wedge a_1 = b \wedge b_1 = 0$. Now let $x \in \bot(a_1) \cap \bot(b_1)$. Then, by part (2), $\bot(x) \in \top(\bot(a_1)) \cap \top(\bot(b_1)) = \{1_{L^\perp}\}$. Hence $\bot(x) = 1_{L^\perp} = L$, which implies $x = 0$. Therefore,
\[
\bot(a_1) \wedge \bot(b_1) = \{0\}.
\]

(4) Using part (2), for every $a, b \in L$, we have $b \in \bot(a)$ if and only if $\bot(b) \in \top(\bot(a))$. Moreover, for every $a, b \in L$, we can write
\[
\bot(a) \cap \bot(b) = \{0\}
\Longleftrightarrow
\top(\bot(a)) \cap \top(\bot(b)) = \top(\bot(a) \cap \bot(b)) = \{0\}^\top = L = 1_{L^\perp}.
\]
Therefore, we are done.
\end{proof}
The skeleton of a pseudocomplemented lattice  $L$, $\text{Skel}L=\{a^*: a\in L\}$, was defined in \cite[Theorem 100]{G0}.  It was proved that  the ordering of $L$ orders $\text{Skel}L$ and makes it as the Boolean lattice as the follows: \[\text{For}~a, b\in L, a^*\wedge b^*=(a\vee b)^*~\text{and}~ a^*\vee_{\text{Skel}L} b^*=(a\wedge b)^*.\]For simplicity, we will use the symbol ``$L^*$'' instead of the symbol ``$\text{Skel}L$'' from now on.
\begin{thm}\label{dube1}
The following statements hold.
\begin{enumerate}
\item If $L$ is a pseudocomplemented lattice, then $L^\perp$ and $L^*$ are two isomorphic lattices and so $L^\perp$ is a Boolean lattice.

\item If $L$ is a Boolean algebra, then $L^\perp$ is too.
\end{enumerate}
\end{thm}
\begin{proof}
(1) Define $\phi: L^\perp\longrightarrow L^*$ by $\phi(\bot(a))=a^*$. We claim that $\phi$ is a lattice isomorphism. First assume $\bot(a)=\bot(b)$ for some $a, b\in L$. Then, clearly, $a^*=\max(\bot(a))=\max(\bot(a))=b^*$. This shows that $\phi$ is a well-defined map. Trivially, $\phi$ is an onto map. To see $\phi$ is one-one, assume $\phi(\bot(a))=\phi(\bot(b))$ for some $a, b\in L$. Then $a^*=b^*$. This implies  $\bot(a)={\dto}a^*={\dto}b^*=\bot(b)$. It is enough to show $\phi$ preserves finite joint  and finite meet.   Consider $\bot(a)$ and $\bot(b)$ in $L^\perp$. Then, we have \[\phi(\bot(a)\wedge \bot(b))=\phi(\bot(a\vee b))=(a\vee b)^*=a^*\wedge b^*=\phi(\bot(a))\wedge\phi(\bot(b)),~\text{and}\]\[\phi(\bot(a)\vee_{L^\perp} \bot(b))=\phi(\bot(a\wedge b))=(a\wedge b)^*=a^*\vee_{L^*} b^*=\phi(\bot(a))\vee_{L^*}\phi(\bot(b)).\] Therefore, two lattices are isomorphic and hence $L^\perp$ is a Boolean algebra.

(2) This follows from (1).
\end{proof}
The converse of the theorem above Part (2) need not be true. It is easy to see that a lattice $L$ is a dense lattice \ifif $L^\perp=\{0, 1\}$. Now, let $L$ be a dense lattice which is not a Boolean algebra (e.g., Example \ref{21}(2)). Then, $L^\perp$ is a Boolean algebra. 
\begin{Rem}\label{dast}
Let $L_1$ and $L_2$ be lattices. Recall that a function $f: L_1 \to L_2$ is said to be a \emph{dual homomorphism} if for every $a, b \in L_1$,
\[
f(a \vee b) = f(a) \wedge f(b)
\quad \text{and} \quad
f(a \wedge b) = f(a) \vee f(b).
\]
If $f: L_1 \to L_2$ is a lattice epimorphism (respectively, a dual lattice epimorphism), then $\theta_f$ is a congruence on $L_1$, and the map $g: L_1 / \theta_f \to L_2$ defined by $g([a]_{\theta_f}) = f(a)$ is an isomorphism (respectively, a dual isomorphism).

Now consider the map $\phi: L \to L^\perp$ defined by $\phi(a) = \bot(a)$, and the map $\varphi: L \to L^\top$ defined by $\varphi(a) = \top(a)$. It is easy to see that $\phi$ is a dual lattice epimorphism, while $\varphi$ is a lattice epimorphism. Therefore,
\[
L^\perp \cong L^{\mathrm{op}} / \theta_\phi
\quad \text{and} \quad
L^\top \cong L / \theta_\varphi.
\]
\end{Rem}

\subsection{The lattice of annihilator ideals of a distributive lattice}
We denote by $\mathcal{A}(L)$ the set of all annihilator ideals of a lattice $L$. It is well-known that $\mathcal{A}(L)=\{\bot(I): I\in I(L)\}=\{I^*: I\in I(L)\}=\text\\{Skel}~I(L)$, partially ordered with the inclusion, is a lattice with the following operations:
\[\bot(I)\wedge \bot(J)=\bot(I)\cap \bot(J)=\bot(I\vee J)~\text{and}~ \bot(I)\vee_{\mathcal{A}} \bot(J)=\bot(I\cap J).\] Indeed, we have $\mathcal{A}(L)=(I(L))^\perp$. So, $\mathcal{A}(L)$ is the $\perp$ of the lattice $I(L)$. Hence, by Theorem~\ref{dube1}, $\mathcal{A}(L)$ is a Boolean algebra. The following result follows from Proposition \ref{sadr}.
\begin{lem}
For a distributive lattice $L$ the following statements are equivalent.
\begin{enumerate}
\item $\mathcal{A}(L)$ is a sublattice of the lattice $I(L)$.

\item $I(L)$ is a co-normal lattice.

\item For each  $I, J\in I(L)$, $\bot(\bot(I)\wedge \bot(J))=\perp^2(I)\vee \perp^2(J)$.  
\end{enumerate}
\end{lem}
\begin{prop}\label{mokh}
For a distributive lattice $L$ the following statements hold.
\begin{enumerate}
\item The lattice $L^\perp$ is a sublattice of the lattice $\mathcal{A}(L)$.

\item The lattice $\mathcal{A}(L)$ is a quotient of the lattice $I(L)^{op}$. 
\end{enumerate}
\end{prop}
\begin{proof}
(1) Since for each $a\in L$, we have $\bot(a)=\bot(\downarrow a)$, it follows that $L^\perp\subseteq\mathcal{A}(L)$. Now, assume $a, b\in L$. Then, \[\bot(a)\vee_{\mathcal{A}} \bot(b)=\bot(\downarrow a)\vee_{\mathcal{A}}  \bot(\downarrow b)=\bot(\downarrow a\wedge \downarrow b)=\bot(\downarrow (a\wedge b))=\bot(a\wedge b)\in L^\perp, ~\text{and}\] \[\bot(a)\wedge \bot(b)=\bot(\downarrow a)\wedge \bot(\downarrow b)=\bot(\downarrow a\vee \downarrow b)=\bot(\downarrow (a\vee b))=\bot(a\vee b)\in L^\perp.\]
Thus $L^\perp$ is closed under the operations of $\mathcal{A}(L)$, and hence is a sublattice.

(2) Since  $\mathcal{A}(L)=(I(L))^\perp$, it follows from Remark \ref{dast} that
$\mathcal{A}(L)$ is a quotient of $I(L)^{op}$.
\end{proof}
\begin{lem}\label{nazari}
Let $\phi: L_1\longrightarrow L_2$ be a lattice isomorphism, $a\in L_1$ and $I\in I(L_1)$. Then, the following statements hold.
\begin{enumerate}
\item $\phi(\bot(I))=\bot(\phi(I))$.

\item $\phi(\bot(a))=\bot(\phi(a))$.
\end{enumerate}
\end{lem}
\begin{proof}
The proof is straightforward.
\end{proof}
\begin{thm}
Let $\phi: L_1 \to L_2$ be a lattice isomorphism. Then the following statements hold:
\begin{enumerate}
\item  The map $\phi^\perp: L_1^\perp \to L_2^\perp$ defined by $\phi^\bot(\bot(x)) = \bot(\phi(x))$ is a lattice isomorphism. Hence, $L_1^\perp$ and $L_2^\perp$ are isomorphic.

\item $\mathcal{A}(L_1)$ and $\mathcal{A}(L_2)$ are isomorphic.

\item If $\perp_{L_1}: L_1 \to L_1^\perp$ and $\perp_{L_2}: L_2 \to L_2^\perp$ are defined by $\perp_{L_1}(x) = \bot(x)$ and $\perp_{L_2}(x) = \bot(x)$, respectively, then the following diagram is commutative.
\[
\begin{tikzcd}
L_1 \arrow[r, "\phi"] \arrow[d, "\perp_{L_1}"'] & L_2 \arrow[d, "\perp_{L_2}"] \\
L_1^\perp \arrow[r, "\phi^\perp"'] & L_2^\perp
\end{tikzcd}
\]
\end{enumerate}
\end{thm}
\begin{proof}
(1) Let $\bot(x), \bot(y) \in L_1^\perp$ with $\bot(x) = \bot(y)$. By Lemma~\ref{nazari}, $\bot(x) = \bot(y)$ implies $\bot(\phi(x)) = \bot(\phi(y))$. Hence $\phi^\bot(\bot(x)) = \phi^\bot(\bot(y))$, so $\phi^\perp$ is well-defined. We claim that $\phi^\perp$ is a lattice isomorphism.

Since $\phi$ is injective, it follows that $\phi^\perp$ is injective. To see that $\phi^\perp$ is surjective, let $\bot(b) \in L_2^\perp$ be arbitrary. Since $\phi$ is surjective, there exists $a \in L_1$ such that $\phi(a) = b$. Then $\phi^\bot(\bot(a)) = \bot(\phi(a)) = \bot(b)$. Thus $\phi^\perp$ is surjective.

Now let $\bot(a), \bot(b) \in L_1^\perp$. By Lemma~\ref{nazari}, we have
\[
\phi^\bot(\bot(a) \wedge \bot(b))
= \phi^\bot(\bot(a \vee b))
= \bot(\phi(a \vee b))
= \bot(\phi(a) \vee \phi(b))
= \bot(\phi(a)) \wedge \bot(\phi(b))
= \phi^\bot(\bot(a)) \wedge \phi^\bot(\bot(b)),
\]
and
\[
\phi^\bot(\bot(a) \vee \bot(b))
= \phi^\bot(\bot(a \wedge b))
= \bot(\phi(a \wedge b))
=\bot(\phi(a) \wedge \phi(b))
= \bot(\phi(a)) \vee \bot(\phi(b))
= \phi^\bot(\bot(a)) \vee \phi^\bot(\bot(b)).
\]
Therefore $\phi^\perp$ preserves both operations, and hence is a lattice isomorphism.

(2) Since $L_1$ and $L_2$ are isomorphic, their ideal lattices $I(L_1)$ and $I(L_2)$ are also isomorphic. Applying part (1) to these isomorphic lattices, we obtain that
\[
\mathcal{A}(L_1) = (I(L_1))^\perp
\quad \text{and} \quad
\mathcal{A}(L_2) = (I(L_2))^\perp
\]
are isomorphic.

(3) It suffices to show that $\perp_{L_2} \circ \phi = \phi^\perp \circ \perp_{L_1}$. Let $a \in L_1$. Then
\[
(\perp_{L_2} \circ \phi)(a)
= \perp_{L_2}(\phi(a))
= \bot(\phi(a))
= \phi^\bot(\bot(a))
= \phi^\bot(\perp_{L_1}(a))
= (\phi^\perp \circ \perp_{L_1})(a).
\]
Hence the diagram commutes.
\end{proof}
Let $X$ and $Y$ be homeomorphic spaces. Then $C(X)$ and $C(Y)$ are
isomorphic rings. Furthermore, by \cite[Theorem 4.4]{BB1} and
\cite[Theorem 4.3]{TA}, the lattices $Z(X)$ and $Z(Y)$, as well as
$Z^\circ(X)$ and $Z^\circ(Y)$, are isomorphic. Consequently, by the
preceding theorem, we obtain the following result. 
\begin{cor}
Let $X$ and $Y$ be two homeomorphic spaces. Then,
\begin{enumerate}
\item  Two lattices $Z(X)^\perp$ (resp., $Z^\circ(X)^\perp$) and $Z(Y)^\perp$ (resp., $Z^\circ(X)^\perp$) are isomorphic.

\item Two lattices $\mathcal{A}(Z(X))$ (resp., $\mathcal{A}(Z^\circ(X))$) and $\mathcal{A}(Z(Y))$ (resp., $\mathcal{A}(Z^\circ (Y)$) are isomorphic.
\end{enumerate}
\end{cor}
This observation might suggest that a lattice $L$ has property $\mathcal{P}$ if and only if $L^\perp$ has the dual counterpart of $\mathcal{P}$. Nevertheless, the following example shows that this is not true in general. 
\begin{exam}
It is well-known that $X$ is a normal space \ifif $\mathcal{C}(X)$ is a co-normal lattice, see \cite[Proposition 4.4]{C}. Also, by \cite[Corollary 4.1.3]{ST}, $X$ is a normal lattice \ifif $I(\mathcal{C}(X))$ is a normal lattice.  Now, let $X$ be a non-normal space (e.g., the Moore plane, see \cite[15.3.b]{W}). Then, the lattice $L=\mathcal{C}(X)$ and also the lattice $I(L)$ are not co-normal lattices. However, by Proposition \ref{mokh}, $\mathcal{A}(L)=I(L)^\perp$ is a Boolean algebra and hence, it is a normal lattice. This shows that the reverse of the Theorem \ref{dube} need not be true. 
\end{exam}

\end{document}

\begin{defn}
Let $L$ be a sublattice of a Boolean algebra $B$. We say that $L$ is a $\emph{semi-$\vee$-lattice}$ in $B$ if, for every nonzero elements $a,b\in L$ with $b\not\leq a$, there exists a nonzero element $x\in L$ such that
$
x\leq a^\prime\wedge b.
$
\end{defn}

It is easily seen that every $\vee$-dense lattice in a Boolean algebra is semi-$\vee$-dense in that Boolean lattice. However, the converse does not hold in general, as shown by the following lemma.
\begin{lem}\label{15}
Assume that $X$ is a Tychonoff space. 
\begin{enumerate}

\item  $Z(X)$ is a semi-$\vee$-dense in the Boolean algebra $P(X)$ (i. e., the power set of $X$);

\item  $Z(X)$ is a $\vee$-dense in the Boolean algebra $P(X)$ if and only if every point of $X$ is a $G_{\delta}$-point.
\end{enumerate}
\end{lem}
\begin{proof}
(1)  It is clear.

(2)$\Rightarrow$ Suppose that $x\in X$, then $\{x\}=\bigvee_{i\in I}Z(f_i)$ and so $\{x\}=Z(f_j)$ for some $j\in I$.

(2)$\Leftarrow$ Let $A\in P(X)$. Clearly, $A=\bigcup_{a\in A}\{a\}$. Since $\{a\}\in Z(X)$, for every $a\in A$, we are done.
\end{proof}

\begin{prop}\label{16}
Let $L$ be semi-$\vee$-dense in a Boolean algebra $B$. Then, the following statements hold:
\begin{enumerate}
\item For each $a\in L$, $\downarrow a = a^{\perp\perp}$.
\item There exists $b\in a^\perp$ such that $b^\perp = a^{\perp\perp}$ if and only if $a$ is complemented in $L$.
\item The lattice $L$ is a $\perp$-lattice if and only if $L$ is a Boolean algebra.
\end{enumerate}
\end{prop}
(1) The inclusion $\downarrow a \subseteq a^{\perp\perp}$ is trivial. Conversely, let $b\in a^{\perp\perp}$. If $b\not\leq a$, then in the Boolean algebra $B$, we have $b\wedge a^c \neq 0$. By semi-$\vee$-denseness, there exists a non-zero $x\in L$ with $x\leq b\wedge a^c$. Then $x\in a^\perp$ and $x\wedge b\neq 0$, contradicting $b\in a^{\perp\perp}$. Thus $b\leq a$.

(2) ($\Leftarrow$) Suppose $a$ is complemented in $L$, so $a^\prime \in L$. Let $b=a^\prime$. We show $b^\perp=a^{\perp\perp}$. Clearly $b^\perp\supseteq a^{\perp\perp}$. For the reverse inclusion, take $x\in b^\perp$ and $y\in a^\perp$. Since $B$ is Boolean, $y\leq a^\prime=b$. Hence $x\wedge y \leq x\wedge b=0$, so $x\in a^{\perp\perp}$.

($\Rightarrow$) Assume there exists $b\in a^\perp$ with $b^\perp=a^{\perp\perp}$. Since $b\in a^\perp$, we have $b\leq a^\prime$.  By semi-$\vee$-denseness, there is a  non-zero $x\in L$ with $x\leq a^\prime\wedge b^\prime$. Then $x\in b^\perp$ (as $x\wedge b=0$). Also, since $x\leq a^\prime$, we have $x\wedge a=0$, so $x\in a^\perp$. But then $x\notin a^{\perp\perp}$ because $x\wedge x = x \neq 0$ (or more directly: $x\not\perp x$). This contradicts $b^\perp=a^{\perp\perp}$. Therefore $b=a^\prime$, so $a$ is complemented in $L$.

(3) This follows immediately from part (2).

The following well-known lemma will be needed in the sequel.
\begin{lem}\label{00}
Let $U$ be an open subset of $Y$ and $X$ be a dense subspace of $Y$. Then, $\cl_{Y}(U\cap X)=\cl_{Y}U$. 
\end{lem}

\section{A new class of topological spaces}
\begin{defn}
A topological space $X$ is called an $\emph{DIC-space}$ if for every pair of disjoint open sets $U, V\subseteq X$, there exists two closed subsets $F, G$ of $X$ such that $U\subseteq F$, $V\subseteq G$ and $\inte_X F\cap\inte_X G=\emptyset$.
\end{defn}
The following easy lemma give a more clear of a $DIC$-space.
\begin{lem}\label{0}
The following statements hold.
\begin{enumerate}
\item A space $X$ is $DIC$.

\item Every two disjoint open sets can be separated by two disjoint regular open sets.

\item For each pair of disjoint open subsets $U, V$ of $X$, we have $\inte\cl U\cap\inte\cl V=\emptyset$.
\end{enumerate}
\end{lem}
\begin{proof}
(1)$\Rightarrow$(2) Let $U, V$ be two disjoint open sets in $X$. By hypothesis, there are two closed subsets $A, B$ of $X$ such that $U\subseteq A$ and $V\subseteq B$ and $\inte A\cap\inte B=\emptyset$. Put $G=\inte A$ and $H=\inte B$. Then $G, H$ are two regular open sets and $U\subseteq G$ and $V\subseteq H$. 

(1)$\Rightarrow$(2)
Let $U, V$ be two disjoint open sets in $X$. By hypothesis, there are two regular open sets $G, H$ such that $U\subseteq G$ and $V\subseteq H$.
This implies $\inte\cl U\subseteq G$ and $\inte\cl V\subseteq \inte H$. Hence, $\inte\cl U\cap\inte\cl V=\emptyset$. 

(3)$\Rightarrow$(1) Let $U, V$ be two disjoint open sets in $X$. Then,  $U\subseteq\cl U$, $V\subseteq \cl V$ and by hypothesis,  we have, $\inte\cl U\cap\inte\cl V=\emptyset$. Thus, $X$ is a $DIC$-space.
\end{proof}
Recall from \cite{AT} and \cite{DT} that a space $X$ is a $\emph{WED-space}$ if every pair of disjoint open sets can be separated by two zero-sets whose  interior are disjoint. Trivially, every $WED$-space is an $DIC$-space, hence every extremally disconnected and every perfectly normal space is a $DIC$-space. We see an example of an $DIC$-space which is not a $WED$-space.
\begin{exam}
 Let $D$ be an uncountable discrete space and let $X=D\cup\{\si\}$ be the one-point compactification of $D$.  It was shown in \cite[Example 3.9]{AT} that $X$ is not a $WCF$-space and hence, it is not a $WED$-space. Let $U$ and $V$ be two disjoint open sets in $X$. If $U$ or $V$ is finite, then we are done. Now, let both $U$ and $V$ be infinite and $\si\not\in U$ and $\si\not\in V$. Then, it is easy to see that $\si\in\cl U\cap\cl V$. On the other hand,  $\si\not\in\inte\cl U\cap\inte \cl V$. Since, $\si\in\inte\cl U\cap\inte \cl V$ implies $\si\in G\subseteq\cl U=U\cup\{\si\}$, where $G$ is open in $X$. Hence, we must have $D\setminus U\subseteq D\setminus G$. This implies $D\setminus U$ is finite and hence $U\cap \Emptyset\neq V$, a contradiction. Thus we must have, $\inte\cl U\cap\inte \cl V=\emptyset$. Finally, we assume both $U$ and $V$ are infinite and $\si\in U$ ($\si\in V$). Then, $V\cap U=\emptyset$ implies $V\subseteq D\setminus U$. But, $D\setminus U$ is finite (since $\si\in U$ and $U$ is open). Thus, $V$ is finite, a contradiction.  Similarly, if $\si\in V$, we obtain $U$ is finite,   a contradiction. Therefore, by Lemma \ref{0}, $X$ is a $DIC$-space.
\end{exam}
\begin{thm}
A space $X$ is a $DIC$-space \ifif for every pair of ideals $I, J$ of $C(X)$ with $I\cap J=0$, there are two ideals $I_1, J_1$ such that $I_1+ J_1$ is essential and $I\cap I_1=J\cap J_1=0$.
\end{thm}
\begin{proof} $\Rightarrow$ Let $I$ and $J$ be two ideals of $R$ with $I\cap J=0$. Put $A=\bigcup\Coz I$ and $B=\bigcup\Coz J$. Then, $A, B$ are two open sets in $X$ and $I\cap J=0$ implies $A\cap B=\emptyset$. By hypothesis, there are two closed sets $G, H$ in $X$ such that $A\subseteq G$,  $B\subseteq H$ and $\inte_X G\cap\inte_X H=\emptyset$. By completely regularity of $X$, there are two subsets $S, K$ of $C(X)$ such that $G=\bigcap_{f\in S}Z(f)$ and $H=\bigcap_{f\in K}Z(f)$. Set, $I_1=<S>, J_1=<K>$. Then, we have \[\inte_X(\bigcap Z[I_1+ J_1])=\inte_X(\bigcap Z[I_1])\cap\inte_X(\bigcap Z[J_1])=\inte_X(\bigcap Z[S])\cap\inte_X(\bigcap Z[K])=\inte_X G\cap\inte_X H=\emptyset.\] This implies $I_1 +J_1$ is an essential ideal of $C(X)$. Moreover, $A\subseteq G$ and  $B\subseteq H$ imply $\bigcup\Coz I\subseteq\bigcap Z[I_1]$ and $\bigcup\Coz J\subseteq\bigcap Z[J_1]$, respectively. Thus, we have \[\bigcap Z[I_1\cap I]=\bigcap Z[I]\cup \bigcap Z[I_1]=X~\text{and}~\bigcap Z[J_1\cap I]=\bigcap Z[I]\cup \bigcap Z[J_1]=X.\] Hence, $I\cap I_1=J\cap J_1=0$.

$\Leftarrow$ Let $U$ and $V$ be two disjoint open sets in $X$. By completely regularity of $X$, there are  two ideals $I$ and $J$ of $C(X)$ such that $U=\bigcup\Coz I$ and $B=\bigcup\Coz J$. Since $U\cap V=\emptyset$, it follows  $I\cap J=0$. By hypothesis, there are two ideals $I_1, J_1$ of $C(X)$ with $I_1+ J_1$ is essential and $I\cap I_1=J\cap J_1=0$. Set $G=\bigcap Z[I_1]$ and $H=\bigcap Z[J_1]$. Then, from $I\cap I_1=J\cap J_1=0$, we have $U=\bigcup\Coz I\subseteq\bigcap Z[I_1]=G$ and $V=\bigcup\Coz J\subseteq\bigcap Z[J_1]=H$. Moreover, the essentiality of  $I_1+ J_1$ implies that \[\inte_X G\cap\inte_X H=\inte_X(\bigcap Z[I_1])\cap \inte_X(\bigcap Z[J_1])=\inte_X\bigcap Z[I_1+ J_1]=\emptyset.\]
\end{proof}
\begin{thm}
The following statements hold.
\begin{enumerate}
\item Let $X$ be a dense subspace of the space $Y$. Then, $X$ is a $DIC$ \ifif $Y$ is a $DIC$-space.

\item Let $X$ be a completely regular Hausdorff space. Then, $X$ is a $DIC$-space \ifif $\beta X$ is a $DIC$-space.

\item Let $X\subseteq Y\subseteq\beta X$. Then $X$ is a $DIC$-space \ifif $Y$ is a $DIC$-space.

\item A  space $X$ is a $DIC$-space \ifif $\upsilon X$ is a $DIC$-space.
\end{enumerate}
\end{thm}
\begin{proof}
(1) $\Rightarrow$ Let $U$ and $V$ be two disjoint open sets in $Y$. Then $U\cap X\cap V\cap X=\emptyset$. By hypothesis, there are two closed subsets $A, B$ in $X$ such that $U\cap X\subseteq A$, $V\cap X\subseteq B$ and $\inte_X A\cap\inte_X B=\emptyset$. There are two closed subsets $A_1, B_1$ in $Y$ such that $A=A_1\cap X$ and $B=B_1\cap X$. By Lemm \ref{00}, we have \[U\subseteq \cl_{Y}U=\cl_{Y}(U\cap X)\subseteq\cl_{Y}A\subseteq A_1~\text{and}~V\subseteq \cl_{Y}V=\cl_{Y}(V\cap X)\subseteq\cl_{Y}B\subseteq B_1.\] We claim that $\inte_{Y} A_1\cap\inte_{Y} B_1=\emptyset$. To see it, by the contrary, let $\inte_{Y} A_1\cap\inte_{Y} B_1\neq\emptyset$. Then, $\inte_{Y} A_1\cap\inte_{Y} B_1\cap X\neq\emptyset$, since $X$ is dense in $Y$. Hence, there exists $p\in X$ with $p\in \inte_{Y} A_1\cap\inte_{Y} B_1$. Thus,  there are open sets $H, K$ in $Y$ such that $p\in H\subseteq A_1$ and $p\in K\subseteq B_1$. These imply $p\in H\cap X\subseteq A_1\cap X=A$ and $p\in K\cap X\subseteq B_1\cap X=B$. Thus, $p\in \inte_X A$ and $p\in\inte_X B$, i.e., $p\in\inte_X A\cap\inte_X B$, a contradiction. Therefore, $Y$ is a $DIC$-space. 

$\Leftarrow$ Let $U$ and $V$ be two disjoint open sets in $X$. Then, there are two disjoint open sets $U_1$ and $V_1$ in $Y$, respectively such that $U=U_1\cap X$ and $V=V_1\cap X$, since $X$ is dense in $Y$. By hypothesis, there are two closed subsets $A, B$ in $Y$ such that $U_1\subseteq A$, $V_1\subseteq B$ and $\inte_{Y}A\cap\inte_{Y}B=\emptyset$. This implies $U\subseteq A\cap X$ and $V\subseteq B\cap X$. We claim that $\inte_X(A\cap X)\cap\inte_X(B\cap X)=\emptyset$. Assume $x\in \inte_X(A\cap X)\cap\inte_X(B\cap X)$. Then, there are two open subsets $G, H$ in $X$ such that $x\in G\subseteq A\cap X$ and $x\in H\subseteq B\cap X$. There are open subsets $G_1, H_1$ in $Y$ such that $G=G_1\cap X$ and $H=H_1\cap X$. Hence, we have \[x\in G_1\subseteq \cl_{Y}G_1=\cl_{Y}(G_1\cap X)\subseteq A~\text{and}~x\in H_1\subseteq\cl_{Y}H_1=\cl_{Y}(H_1\cap X)\subseteq B.\] That is $x\in \inte_{Y}A\cap\inte_{Y}B$, a contradiction. Therefore, $X$ is $DIC$-space.

(2) This follows from (1). 

(3) By \cite[Theorem 6.7]{GJ}, $\beta Y=\beta X$. Thus, by Part (2), $X$ is a $DIC$-space \ifif $\beta Y$ is so, and so again by Part (2), $X$ is a $DIC$-space \ifif $Y$ is so.

(4) This follows from Part (3).
\end{proof}